\documentclass[12pt,reqno]{amsart}
\usepackage{latexsym, mathrsfs} 
\usepackage{calrsfs}
\usepackage{psfrag}
\usepackage{amsfonts}
\usepackage{amssymb,stmaryrd}
\usepackage{wrapfig, graphicx}
\usepackage{amscd}
\usepackage{epsfig}
\newtheorem{theorem}{Theorem}
\newtheorem{lemma}[theorem]{Lemma}

\newtheorem{corollary}[theorem]{Corollary}
\newtheorem{proposition}[theorem]{Proposition}
\newtheorem{defin}[theorem]{Definition}

\newenvironment{proofof}[1]{\noindent {\bf Proof of #1.}}{ \hfill\qed\\ 
}
\renewenvironment{itemize}{\begin{list}{$\bullet$} 
{\setlength{\labelwidth}{1cm}
\setlength{\leftmargin}{0.4cm}
\setlength{\rightmargin}{0.3cm}
\setlength{\itemindent}{0mm}
\setlength{\parsep}{0.2ex}
\setlength{\topsep}{0.5ex}
\setlength{\itemsep}{0ex}
\slshape}}{\end{list}}

\newlength{\figboxwidth}
\newcommand{\abx}{(X,\omega)}
\newcommand{\aby}{(Y,\tau)}

\newcommand{\slr}{ {\rm SL}_2(\R) }
\newcommand{\glr}{ {\rm GL}_2(\R) }

\newcommand{\slz}{ {\rm SL}_2(\Z) }

\newcommand{\slv}[1]{ {\rm SL}(#1) }

\newcommand{\R}{\mathbb{R}}
\newcommand{\C}{\mathbb{C}}
\newcommand{\Q}{\mathbb{Q}}
\newcommand{\Z}{\mathbb{Z}}
\def\L{\mathbb{L}}
\newcommand{\N}{\mathbb{N}}

\renewcommand{\S}{\mathbb{S}}
\newcommand{\T}{\mathbb{T}}

\newcommand{\sE}{\mathscr{E}}
\newcommand{\sF}{\mathscr{F}}

\newcommand{\sO}{\mathscr{O}}

\newcommand{\CB}{\mathcal{B}}
\newcommand{\CC}{\mathcal{C}}

\newcommand{\CF}{\mathcal{F}}

\newcommand{\CH}{\mathcal{H}}

\newcommand{\CL}{\mathcal{L}}

\newcommand{\CO}{\mathcal{O}}

\newcommand{\CT}{\mathcal{T}}
\newcommand{\CU}{\mathcal{U}}

\def\e{\varepsilon}

\def\x{\bar{x}}
\def\y{\bar{y}}
\def\Rec{\mathcal{R}}
\def\diact{\circ_{\scriptscriptstyle{d}}}

\renewcommand{\d}{{\rm d}}

\DeclareMathOperator{\1}{1}

\DeclareMathOperator{\aff}{Aff^{+}}
\DeclareMathOperator{\aut}{Aut}
\DeclareMathOperator{\vol}{vol}

\DeclareMathOperator{\D}{D}

\DeclareMathOperator{\id}{id}
\DeclareMathOperator{\lcm}{lcm}
\DeclareMathOperator{\per}{Per}

\DeclareMathOperator{\pr}{pr}

\DeclareMathOperator{\im}{Im}
\DeclareMathOperator{\re}{Re}

\title[ ]{Modular fibers and illumination problems.}  
\author{P.~Hubert, M.~Schmoll and S.~Troubetzkoy}
\address{Laboratoire d'Analyse, Topologie et Probabilités\\
Fédération de Recherches des Unités de Mathématiques de Marseille and\\
Université Paul Cézanne\\
Case Cour A, Avenue Escadrille Normandie-Niemen, F-13397 Marseille Cedex 20, France}

\email{hubert@cmi.univ-mrs.fr}

\address{ University of Notre Dame, Department of Mathematics,
255 Hurley Hall, Notre Dame, IN 46556-4618
}

\email{schmoll.2@nd.edu}

\address{Centre de physique théorique\\
Institut de mathématiques de Luminy\\
Fédération de Recherches des Unités de Mathématiques de Marseille and\\
Université de la Méditerranée\\
Luminy, Case 907, F-13288 Marseille Cedex 9, France}
\email{troubetz@iml.univ-mrs.fr}
\urladdr{http://iml.univ-mrs.fr/{\lower.7ex\hbox{\~{}}}troubetz/}

\date{\today}
\begin{document}

\begin{abstract}
For a Veech surface $\abx$, we characterize $\aff \abx$ invariant 
subspaces of $X^n$ and prove that non-arithmetic 
Veech surfaces have only finitely many invariant subspaces 
of very particular shape (in any dimension). 
Among other consequences we find copies of $\abx$ embedded 
in the moduli-space of translation surfaces. 
We study illumination problems in (pre-)lattice surfaces.
For $(X,\omega)$ prelattice we prove the at most countableness of
points
non-illuminable from any $x \in X$.  Applying our results on invariant
subspaces  we prove the
finiteness of these sets  when $\abx$ is Veech.
\end{abstract}
\maketitle
\section{Introduction}

Consider a room (a plane domain) with mirror walls and a point
source of light which emits rays in all directions.  Is the
whole room illuminated?  This question is known as 
{\it the illumination problem}. 
The problem is attributed to Ernst Strauss in the 
1950's.  It 
has appeared in various lists of unsolved problems
\cite{Kl2,KlWa,CrFaGu}.  Certain publish versions of the
question specify that the room has polygonal
boundary \cite{Kl1,Kl2}, however the earliest published version
of the question and results which we found is for
rooms with smooth boundary \cite{Pe}. 

In this article we are interested in this question when the
room is a polygon.  So far, except for the trivial fact that convex rooms 
are illumined by any  point, all the known results on
the illumination problem in polygons are negative.  Tokarsky
has constructed  polygons $P$ and point sources $p \in P$ which
do not illumine all points $q \in P$ \cite{To}.\footnote{A similar
construction was already evident in an earlier unpublished letter of
M.~Boshernitzan to H.~Masur \cite{Bo}.}
There have been various generalizations of the illumination
property, for example the study of illumination by search lights
\cite{ChGa} and the study of the finite blocking property
\cite{Mo1,Mo2,Gu,HiSn}.

Consider a rational polygon $P$.  There is a well known procedure of
unfolding $P$ to a flat surface with conic singularities.  
We will state our results in terms of certain
classes of flat surfaces with singularities.  By this unfolding
procedure they hold for the polygons which unfold to them.  In 
particular we consider two classes of surfaces, Veech surfaces and
more generally prelattice surfaces.  Their definition 
will be given below.  
Here we just note that all regular polygons
unfold to Veech surfaces.

For surfaces we say that $p \in X$ illumines $q \in X$ if there is a
geodesic connecting $p$ to $q$. Our main illumination results are

\begin{theorem}\label{thm1}
Let $X$ be a prelattice surface.  Then for any points $p \in X$,
the  set of points
$q \in X$ which are not illuminable from  $p$
is at most countable.
\end{theorem}
To prove Theorem \ref{thm1} we establish a quantitative version of Kronecker's theorem.
\begin{theorem}\label{Schmoll}
Let $X$ be a Veech surface.  Then for any points $p \in X$,
the  set of points $q \in X$ which are not illuminable from  $p$
is finite. 
\end{theorem}
\noindent
Affine homeomorphisms are real linear, i.e.\ they 
preserve geodesic segments on $X$,  hence preserve 
illuminable configurations. Any two points 
$p,q \in X$ located in a convex, open set    
$\CO \subset X$ can be connected by a straight line, 
consequently $p$ and $q$ illuminate one another if 
they do not belong to the {\em exceptional set}
 \begin{equation}\sE:= \big \{(p,q) \in X \times X: 
(\aff (X) \cdot p, \aff(X) \cdot q) \not \in \CO \times \CO  \big \},   
\end{equation} 
which is closed and $\aff(X)$ invariant. Here $\aff(X)$ is the group of affine homeomorphisms of $X$.  Theorem \ref{Schmoll} follows by combining Theorem \ref{thm1} with a description of 
$\aff(X)$-invariant subspaces in $X^2:= X\times X$ for a Veech surface $X$. 
\medskip\\
Our results on invariant subspaces of $X^2$ are of independent 
interest and will cover the first part of the paper. To state the 
Theorems we review the basic notions of translation surfaces, details can
be found in the references \cite{MaTa} and \cite{zo1}.
\medskip\\
A {\em translation surface} is a compact orientable surface with an atlas 
such that away from finitely many points called {\em singularities} all 
transition functions are translations.  
Each holomorphic 1-form $\omega$ on a Riemann surface $X$ induces a
translation structure on the surface by taking natural charts: 
\[ \int^{z}_{z_0}\omega \quad z_0,z \in X\backslash Z(\omega).\] 
The set $Z(\omega)$ where $\omega$ {\em vanishes} represents 
the singularities or {\em cone points} of the translation structure. 
An {\em Abelian differential} 
$\abx$ is a pair consisting of a Riemann surface $X$ and a 
holomorphic $1$-form $\omega$ on $X$.
The notion of  translation surface and Abelian differential 
are equivalent.

The group $\slr$ acts on a translation surface $\abx$ in the following sense.
Given an element $A \in \slr$, we can postcompose
the coordinate functions of the charts of the (translation) atlas of 
$\abx$ by $A$.  It is easy to see that this again yields a 
translation surface, denoted by $A \cdot \abx$.
\medskip\\
We say that $\phi$ is an {\it affine homeomorphism} 
of a translation surface $X$ is $\phi$
is a (orientation preserving) homeomorphism that is a 
diffeomorphism on $X \backslash Z(\omega)$, 
whose differential is a constant element of $\slr$ in each chart of the atlas.
The orientation preserving affine homeomorphisms of a translation 
surface form a group $\aff \abx$.
This group acts naturally on the surface $X$ and diagonally on 
$X^n$: for $\phi \in \aff \abx$
\[ X^n \ni (z_1,..,z_n) \longmapsto \phi \diact (z_1,...,z_n)
:=(\phi(z_1),...,\phi(z_n)) \in X^n.
\]
In this article we will mainly study the 
{\em diagonal action} of $\aff \abx$. 

Our starting point is the so called
{\em Ratner Theorem} of Eskin, Marklof and Morris-Witte\\
{\bf Theorem} \cite{EsMaWi}  
{\em 
If $\abx$ is a Veech surface and $(z_1,\dots,z_n) \in X^n$, 
then the closures 
$\overline{\slv{X,\omega}\diact (z_1,\dots,z_n)}\subset X^n$ are  
linear, complex-algebraic spaces.}
\medskip

\noindent Let $\D: \aff \abx \to \slr$ assign to each $\phi \in \aff \abx$ 
its (constant) differential.
The map $\D$ is a differential homomorphism whose range 
$\slv{X,\omega}$ is a Fuchsian group called the {\it Veech group} 
\cite{v1} of $X$. 
Given an Abelian differential $\abx$, there is the exact sequence
\begin{equation}\label{exactaff} 
\begin{CD}
1 @>>> \aut \abx  @>{i}>> \aff \abx @>{\D}>> \slv{X,\omega} @>>> 1.
\end{CD}
\end{equation}  
The kernel $\aut \abx =\ker (\D)$ is the group of 
{\em automorphisms} of $\abx$, i.e. the group of 
biholomorphic maps of $X$ preserving $\omega$. 
In particular: $\aut \abx$ is finite if $g(X)>1$, 
or if $\abx$ is a torus marked in (at least) one point.     
Another consequence of Equation \eqref{exactaff} is that $\slv{X,\omega}$ is the 
stabilizer of the $\slr$ orbit of $\abx$ in the moduli space 
of translation surfaces.  

We call a translation surface $\abx$ {\em Veech- or lattice surface}, 
if $\slv{X,\omega}$ is a {\em lattice} in $\slr$, 
i.e. $\vol(\slv{X,\omega}\backslash \slr)< \infty$.
\medskip\\
Our first subspace result sharpens the above result  
of Eskin, Marklov and Morris-Witte. 
To begin with take an Abelian differential $\abx$ and  
call a subspace $S \subset X^n$ 
 {\em real-linear}, 
if $S$ is locally, in natural coordinates 
\begin{equation}\label{loccord} (p_1,...,p_n) \mapsto 
(z_1,...,z_n):=(\int^{p_i}_{p_{0,i}}\omega, 
...,\int^{p_n}_{p_{0,n}}\omega) \end{equation}
centered at a point $p_0 =(p_{0,1},...,p_{0,n})\in S
$, 
defined by equations of the shape 
\begin{equation}
a_1z_1+...+a_nz_n=0, \mbox{ with }
(a_1,...,a_n) \in \R^n.
\end{equation} 

\begin{theorem}\label{real}
Let $\abx$ be a lattice surface. Then all     
$\aff \abx$ orbit closures in $X^n$ 
are real-linear, complex-algebraic spaces.
\end{theorem}
A consequence of the real linearity of invariant subspaces is:
\begin{theorem}\label{veech-groups} 
Each connected component $S \subset  X^n$ of a $2$-dimensional, invariant subspace carries a natural structure as a Veech-surface 
$(S,\omega_S)$  with Veech group  $\slv{S,\omega_S}$ commensurable to     
the Veech group $\slv{X,\omega}$. 
\end{theorem}
\noindent Furthermore
\begin{theorem}\label{finite}
Assume $\abx$ is a Veech surface. Then 
we have the following dichotomy:
\begin{itemize}
\item either $\abx$ is arithmetic and 
there are infinitely many $\aff \abx$ invariant subspaces 
of any even dimension in $X^n$,  
\item or $\abx$ is not arithmetic and $X^n$ contains 
finitely many $\aff \abx$ invariant subspaces. 
In local coordinates (like in (\ref{loccord})) 
each invariant hypersurface satisfies an equation
\[ \epsilon_1z_1+...+\epsilon_nz_n =0 \]
for some
$\epsilon_1, \dots, \epsilon_n =-1,0,1.$
Lower-dimensional, invariant subspaces are intersections 
of higher-dimensional ones.
\end{itemize}
\end{theorem}

{\bf Horizontal and vertical subspaces.} 
Given a Riemann surface $X$ call a subspace 
$S \subset X^2$ {\em horizontal} ({\em vertical}), 
if for a finite $F \subset X$
\[ S = \coprod_{x \in F} \{x\}\times X \subset X^2, \quad  
   S = \coprod_{x \in F} X \times \{x\} \subset X^2
\mbox{ respectively}. 
\]  
To specify we sometimes write $X^2$ as $X_h \times X_v$,  
and think of a {\em horizontal} and {\em vertical} component. 
Sometimes we call vertical or horizontal 
subspaces simply {\em parallel} subspaces.   

\noindent{\bf Slope of an invariant subspace in $X^2$.} 
As a further consequence of Theorem \ref{real}, 
a linear, connected $\Gamma \subset \aff \abx$  
invariant subspace $S \subset X^2$ has a well defined 
slope. To start with $S$ is (locally) given by 
an equation 
\[ az_h=bz_v, \text{ with } (z_h,z_v) \in X^2 \text{ and } 
(a,b) \in \R^2.\]  
We define its {\em slope} to be  
\[\frac{a}{b} \in \R \cup \{\infty\}.\]
We can generalize and speak of a foliation $\CF_{\alpha}(X^2)$ 
of $X^2$, locally defined by linear equations with slope 
$\alpha \in \R \cup \{\infty\}$. 
\begin{lemma}\label{welldefined}
The slope is well defined for 
$2$-dimensional $\aff \abx$ invariant 
subspaces of $X^2$.
Moreover the diagonal action of 
$\aff \abx$ preserves any foliation $\CF_{\alpha}(X^2)$.
\end{lemma}

A $2$-dimensional $\aff \abx$ invariant subspace  
of $X^2$ has rational slope (see Theorem \ref{finite} for the non-arithmetic
case).

\noindent Most of our results are consequences of  
Theorem \ref{real} and Theorem \ref{finite}.
\medskip\\

{\bf Primitive and reduced translation surfaces.}
Call an Abelian differential $\abx$ {\em reduced}, if 
$\aff \abx \cong \slv{X,\omega}$. Given a Veech surface $\abx$, 
there is always a reduced Veech surface 
$(X_{red},\omega_{red})$ and a covering map $\pi: X \rightarrow X_{red}$ 
such that 
\[\pi^{\ast} \omega_{red}=\omega \ \text{ and } \ 
\slv{X,\omega} \subset \slv{X_{red},\omega_{red}}.\]
In fact $X^n_{red}$ describes the space of ordered 
$n$-markings on $X$ up to translation maps, 
for more, see the discussion after the introduction.
\medskip\\ 
A {\em primitive} Abelian differential $\abx$ is a translation surface 
which does not admit a translation map to a surface of lower genus, 
i.e. there is no covering map $\pi: \abx \rightarrow (Y,\alpha)$, 
such that $\omega=\pi^{\ast}\alpha$ and $g(X)>g(Y)$.
\medskip\\
{\bf Off-Diagonal subspaces.}  
Given a Riemann surface $S$, we have already 
defined  the {\em diagonal} $D_+$. 
If $\phi: S \rightarrow S$ is an involution, 
we denote its graph by 
\begin{equation}
D_{\phi}:= graph(\phi) :=\{(x,y)\in S \times S: y=\phi(x)\}. 
\end{equation}
If $\abx$ has a unique involution $\phi$, in particular if $\abx$ 
is reduced we simply write $D_{-}$ instead of $D_{\phi}$ and call 
it the {\em off-diagonal}.

The next Theorem classifies all possible invariant subspaces 
under certain assumptions. 
\begin{theorem}\label{classification}
Let $\abx$ be a primitive Veech surface with 
exactly 
\underline{one} cone-point of prime order.  Then 
all connected components of $2$-dimen{\-}sional 
$\slv{X,\omega}$-invariant subspaces $S \subset X^2_{red}$ 
admit a natural translation structure $(S,\omega_S)$, 
with an affine diffeomorphism 
\[ (S,\omega_S) \cong \abx.\]
Furthermore if $\abx$  
is \underline{not} arithmetic, then 
the only $2$-dimensional, linear subspaces in $X^2_{red}$  
away from the parallel ones  are  $D_{+}$ and $D_{-}$. 
The last subspace exists 
if and only if $\abx$ admits an (affine) involution. 
\end{theorem}
\vspace*{1mm}
Since $X^n_{red}$ parametrizes $n$-tuples of marked points 
on $\abx$, it is a {\em fiber}  
\[ X^n_{red} \hookrightarrow \slr \diact \ X^n_{red} 
\stackrel{\pi}{\rightarrow} \slr \cdot \abx \cong 
\slv{X,\omega}\backslash \slr,
\] 
i.e. $\pi^{-1}(\abx) = X^n_{red}$, in the moduli space of $n$-tuple 
markings on deformations of $\abx$, thus we call $X^n_{red}$ the 
{\em modular fiber}. Our results say, 
the modular fiber and hence moduli space itself 
contains Veech surfaces $(S,\omega_S)$ with a 
lattice group $\slv{S,\omega_S} \cong \slv{X,\omega}$. 
In general Veech surfaces $(S,\omega_S) \subset X^2$  
are not isomorphic to the surface $\abx$ defining the 
modular fiber. 
\medskip\\
{\bf More illumination.} 
Given a covering map $\pi: X \rightarrow Y$ we define  
$\pi_2:=\pi \times \pi: X^2 \rightarrow Y^2$ by   
$\pi_{2}(x_1,x_2):=(\pi(x_1),\pi(x_2))$. 
For the following we need the coordinate projection 
$pr_h:X^2 \cong X_h\times X_v \rightarrow X_h$,  
$(p,q) \mapsto p$.
\begin{theorem}\label{nonilluconf}
Take a Veech surface $\abx$ with canonical covering 
$\pi: X \rightarrow X_{red}$. 
Assume that $X^2_{red}$ contains exactly $D_+$ 
and $D_-=D_{\phi}$, $\phi \in \aff(X_{red},\omega_{red})$ 
the (unique) involution, as non-parallel invariant subspaces. 
Then the only non-illuminable configurations  $(p,q)\in X^2$
on $X$ are 
\begin{itemize}
\item all regular pairs in $\pi^{-1}_{2}(D_{-}) $, if 
and only if $pr_h(\pi^{-1}_{2}(D_{+}\cap D_{-})) \subset X$ 
are cone points and    
\item some pairs $(p,q) \in X^2$ of periodic points.  
\end{itemize}
In the later one needs to check for illuminability case by case.
\end{theorem}
\begin{corollary}\label{cornonilluconf}
The previous Theorem applies to all 
known Veech surfaces in genus $2,3$ and $4$ 
with precisely one cone point. 
\end{corollary}
For these examples see C.T.McMullen's papers \cite{mcm1,mcm3}.
\\
{\em Remark}: One can check that 
there are no nonillumination configurations 
in the nonarithmetic genus $2$ (one cone point) Veech case 
by looking at L-shaped representatives and using Möller's 
classification of periodic points.

\section{Background}

For later use we recall some facts on Veech surfaces $\abx$ 
and the flat geometry on $X$ induced by $\omega$.
\medskip\\
A {\em saddle connection} on $\abx$ is a geodesic segment, 
with respect to the flat metric induced by the translation structure, 
starting and ending at zeros of $\omega$. 
A direction on a translation surface is called {\em periodic}
if the flow in this direction is periodic (except for the set
of saddle connections in this direction).  Thus in a periodic direction
the surface decomposes into maximal cylinders bounded by saddle connections.
The {\it modulus} of a cylinder is its width divided by its length.
\medskip\\
An element $A \in \slr \backslash \{\pm 1\}$ is either parabolic, elliptic,
or hyperbolic.  A direction on a translation surface is called {\em  
parabolic}  if there is an affine diffeomorphism that preserves the set
of geodesics in this direction and whose differential is parabolic. 
Veech has shown that a direction is parabolic if and only if 
it is periodic with all moduli commensurable \cite{v1}.  In this case
the action of a parabolic diffeomorphism restricted to a cylinder is
a power of a Dehn twist.
\medskip\\
{\bf Reduced surfaces.} 
As indicated by the exact sequence \eqref{exactaff} of groups, 
we cannot distinguish 
the marked surface $(X,m,\omega)$ 
from the marked surface $(X,\phi(m),\omega)$ for 
any $\phi \in \aut \abx$, at least if we    
consider the {\em moduli space} of marked points 
(equivalence up to translations of $(X,\omega)$).  
The moduli space of marked points is given by the 
quotient $X_{red}:=X/\aut\abx$. 
If $\pi: X \rightarrow X_{red}$ is the quotient map 
we can define the (holomorphic) one form 
\[\omega_{red}:=\pi_{\ast} \omega\]
on $X_{red}$, since all automorphisms preserve $\omega$. 
By definition the Abelian differential $(X_{red},\omega_{red})$  
has no (non-trivial) automorphisms and consequently 
\begin{equation}\label{red} 
\aff(X_{red},\omega_{red})\cong \slv{X_{red},\omega_{red}},  
\end{equation} 
as well as 
\begin{equation}\label{descend} 
\slv{X,\omega} \subset \slv{X_{red},\omega_{red}}  
\end{equation} 
if every affine map of $\abx$ descends to 
an affine map of $(X_{red},\omega_{red})$. 
To see this take  $\psi \in \aff \abx$,  
then $\psi \cdot f \cdot \psi^{-1} \in \aut \abx$ for all 
$f \in \aut \abx$ and consequently 
\[ \pi \circ \psi \circ f(x)=
\pi \circ \psi \circ f \circ \psi^{-1} \circ \psi(x)=\pi \circ \psi(x),
\]
showing that $\psi$ descends. Obviously we have
\begin{equation}\label{finind} 
[\slv{X_{red},\omega_{red}}:\slv{X,\omega}]< \infty.
\end{equation}
A useful observation is 
\begin{proposition}
Given a lattice surface $\abx$ with an involution 
$\phi \in \aff \abx$ and assume $\aff \abx \cong \slv{X,\omega}$. 
Then $\phi$ is the only affine involution of $\abx$. 
\end{proposition}
\begin{proof}
For any affine involution $\phi \in \aff \abx$ we have 
$\D\phi=-\1 \in \slv{X,\omega}$. By assumption   
$\D: \aff \abx \rightarrow \slv{X,\omega}$ is an isomorphism.
\end{proof}
Because the modular fibers $X^n_{red}$ of a lattice surface 
$\abx$ are always reduced and Equation \eqref{descend} holds, 
we can assume all lattice surfaces we are looking at are  
reduced.
\medskip\\

\noindent {\bf Example: Translation surfaces obtained from regular 
$\mathbf{2n}$-gons.}
Veech showed in \cite{v1,v2}, 
that one obtains a lattice surface $\aby$ taking $2$ copies 
of the regular $2n$-gon and identify sides 
of the first regular $2n$-gon with the diametrical, parallel 
sides of the copy. Using just one $2n$-gon and the 
same diametrical gluing scheme, one obtains a  translation 
surface $\abx$ as well. 
Exchanging the two $2n$-gons tiling $\aby$ defines a 
nontrivial automorphism $\phi$ of order $2$ in $\aut \aby$. 
Moreover the quotient of $\aby$ with respect to 
$\aut \aby=<\phi>$ is $\abx$ and thus $(Y_{red}, \tau_{red}) =\abx$.
\medskip\\  
We remark that all presently known 
primitive Veech surfaces of genus $g \geq 2$ 
have one or two cone points, moreover a new 
result of Martin Möller \cite{mm2} shows that the number 
of (algebraically) primitive Veech surfaces with two 
cone points in fixed genus is finite.     
\medskip\\ 
Note, that primitive (lattice) differentials   
are reduced, but a reduced lattice differential is 
not necessary primitive. An infinite set of 
of reduced elliptic differentials $\abx$ with  
$\slv{X, \omega} \cong \slz$ is contained in 
\cite{s2,s3}. Square-tiled surfaces in $\CH(2)$ \cite{hl} 
provide examples of reduced, but not primitive arithmetic 
Veech surfaces too.   
\medskip\\






\section{Invariant subspaces of $X^n$}

The zero dimensional $\slv{X,\omega}$ invariant subspaces of 
$X^n$ are periodic points. Each point $\aby$ 
having finite $\aff \abx$ orbit corresponds to $(X,\omega)$  
marked at some periodic points $(m_1,...,m_n)$.  For 
reduced $\abx$,
\[ [\slv{X,\omega}:\slv{Y,\tau}]= |\slv{X,\omega}\cdot [Y,\tau]|. 
\] 
If $\abx$ is {\em not} an elliptic differential, 
i.e. not a torus cover, a result of Gutkin, Hubert and Schmidt 
\cite{ghs} says there are only {\em finitely many} 
periodic orbits in $\abx$ and thus for all $n$,
there are only finitely many periodic points in $X^n$. 
\medskip\\
{\bf Reality of invariant subspaces and consequences.} 
If a point $z=(z_h,z_v) \in X^2$ is not periodic, at least one 
of its components, say $z_h$, is not periodic and thus has 
closure $X=\overline{\slv{X,\omega}z_h}$ by \cite{ghs}. 
Assume $S$ is a connected $\slv{X,\omega}$ invariant subset 
of $X^2$. Ratner's Theorem  says that 
$\slv{X,\omega}$ invariant subsets in $X^2$ 
are already linear, complex submanifolds. 
We also have
\begin{proposition}
A point $p \in X^2$ is periodic with respect to 
the action of $\slv{X,\omega}$, if and only if 
$p \in L \cap S$ is in the intersection of two 
connected, invariant subspaces $L \neq S \subset X^2$  
with $\dim(L)=\dim(S)=2$. 
\end{proposition}
\begin{proof}
If $L, S \subset X^2$ are two connected, linear, 
invariant subspaces 
of complex dimension $1$, they have 
only finitely many intersection points or $L=S$. By invariance of $L$ and $S$,     
$L \cap S$ is also invariant under the action of $\aff \abx$. 
If $p=(p_h,p_v)\in X \times X$ is periodic, both components 
$p_i$ are periodic and we have 
\[p \in \{p_h\}\times X \cap X \times \{p_v\}.\]
\end{proof}
From now on we call a subspace $S \subset X^n$ 
{\em invariant}, if it is a connected component, 
or a union of connected components 
of an $\slv{X,\omega}$-invariant set.
Let 
\begin{equation}
I=I_m:=\{i_1,...,i_m\} \subset \{1,...,n\}
\end{equation}
be a tuple and 
\[\pr_{I}:X^n \rightarrow X^m=X_{i_1}\times X_{i_2}\times 
...\times X_{i_n}\]  
the canonical projection. Then we have   
\begin{proposition}\label{intersection}
Assume $\abx$ is a lattice surface and $S \subset X^n$ is invariant, 
then
\begin{itemize}
\item intersections of invariant subsets of $X^n$ are invariant and 
\item the sets $\pr_{I_m}(S)$ and 
$\pr^{-1}_{I_m}\circ \pr_{I_m}(S)$ are 
invariant for each index subset $I_{m}$.
\end{itemize}
\end{proposition}
\begin{proof}
The first claim is obvious. For the second claim, we note 
that $\pr_{I_m}(S)$ is $\slv{X,\omega}$ equivariant.
\end{proof}
\vspace*{2mm}

\begin{proof}[{\bf Proof of Theorem \ref{real} for $\mathbf{X^2}$}]
This argument proves the reality of invariant subspaces 
$S \subset X^2$ of dimension $2$. 
We can restrict ourselves to invariant subspaces 
$S$, which are $\aff \abx$ orbit closures, i.e. 
$S$ with a finite number of connected components. 

From \cite{EsMaWi}, see Ratner's Theorem  
we know that any invariant subset $S\subset X^2$ 
which is an $\aff \abx$ orbit-closure of one point $z \in X$ 
is linear, algebraic and of even dimension. In particular 
$S$ has finitely many connected components and is stabilized 
by a finite index subgroup $\Gamma_S \subset \slv{X,\omega}$. 

The statement clearly holds for horizontal 
and vertical subspaces $S$ of $X^2$. Thus we assume  
$S \subset X^2$ is a connected subspace which is 
$\Gamma_S$ invariant for a finite index subgroup 
of $\slv{X,\omega}$ and not horizontal or vertical.  
In this case the two projections 
\[\pi_i:X^2 \rightarrow X_i \cong X, \ i=h,v \] 
are onto and $S$ must be a linear {\em correspondence}, i.e. 
both maps $\pi_i:S\rightarrow X_i$ are finite to one. 
By $\Gamma_S$ invariance and linearity 
we might describe $S$ in a neighborhood of $p \in S$ as $(y,\psi(y))\in X^2$
with an affine map $\psi: \R^2 \rightarrow \R^2$. 
After application of $\phi \in \aff \abx \cap \D^{-1}\Gamma_S$ 
we get by invariance 
\[(\phi(y),\phi \cdot \psi(y))=(z,\psi_{\phi}(z)) \in S,\] 
with an affine map $\psi_{\phi}:\R^2 \rightarrow \R^2$. 
But $z=\phi(y)$ and thus: 
\begin{equation}\label{phipsi} 
(\phi(y),\phi \cdot \psi(y))=(z,\psi_{\phi}(z))=
(\phi(y),\psi_{\phi} \cdot \phi(y)).
\end{equation}
Now for all $\phi \in \aff \abx \cap \D^{-1}\Gamma_S$ 
and all pairs $(\psi,\psi_{\phi})$ 
we obtain by linearity and connectedness of $S$: 
\[\D \psi=\D \psi_{\phi} =: B.\] 
Thus taking the linear part in the second component 
of Equation (\ref{phipsi}) gives  
\[ A \cdot B = B \cdot A \ \mbox{ for all } A \in \Gamma_S.\] 
Without restrictions we assume  
that the vertical and horizontal foliation on $X$ 
is periodic. Taking 
the subgroup of $\Gamma_S$ generated by the Dehn 
twists along the horizontal and vertical direction proves: 
\[ B= b \cdot \id \in \glr 
\ \mbox{ with } b \in \R-\{0\}.\] 
That means all linear, connected $\Gamma_S$ invariant 
$S \subset X^2$ of dimension $2$ are algebraic 
submanifolds of $X^2$ defined by an linear equation 
with real coefficients. 
\end{proof}
\vspace*{1mm}

\noindent Theorem \ref{real} allows to identify directions on horizontal 
and vertical fibers in $X^2$ via any $\slv{X,\omega}$ invariant 
subspace $S$. We can, for example, take a closed 
leaf $\gamma \subset \abx$ viewed as a horizontal embedding 
\[\gamma \subset \abx \rightarrow X \times \{z_v\} \subset X^2.\]
Then an invariant subspace $S$ maps the image of 
$\gamma \subset X \times \{z_v\}$ to a closed leaf 
$\gamma_S \subset \{z_h\} \times X$ {\em in the same direction} 
for any $z_h \in X$. 
\medskip\\
\begin{proofof}{Lemma \ref{welldefined}}
Because $\aff \abx$ acts diagonally on $X^2$ 
and real linear on each component of $X^2$, 
we find for any $\phi \in \aff \abx$ 
\[ a\phi(z_h)-b\phi(z_v)=\D\phi(az_h-bz_v)+c_{\phi}, 
\ \text{with } a,b \in \R, \ c_{\phi} \in \R^2.\]
This implies the slope is well defined for 
$2$-dimensional $\aff \abx$ invariant 
subspaces of $X^2$ and the diagonal action of 
$\aff \abx$ preserves any foliation $\CF_{\alpha}(X^2)$.
\end{proofof}   
\medskip\\
{\bf Canonical holomorphic differential on invariant subspaces.}
A further consequence of the real-linearity of 
an $\slv{X,\omega}$ invariant subspace $S$ is that 
it admits a holomorphic $1$-form $\omega_S$, 
compatible with the $1$-form $\omega$ on $X$ 
with respect to both coordinate projections 
$\pi_i:S\rightarrow X_i$. In fact locally and away 
from its zeros $\omega_S$ might be defined by  
\begin{equation} 
\d z_S:=\alpha\sqrt{1+\alpha^{-2}}\ \d z_h=
\alpha^{-1}\sqrt{1+\alpha^{2}} \ \d z_v
\end{equation}   
for an invariant $S$ with slope $\alpha \neq 0, \infty$. 
If $S$ is horizontal ($\alpha =0$) 
or vertical ($\alpha = \infty$) we have $\d z_S=\d z_h$, 
$\d z_S=\d z_v$ respectively.
\medskip\\
{\bf Cone points of $\mathbf{(S,\omega_S)}$.}
To characterize the cone points of $(S,\omega_S)$ 
and prove finiteness of invariant subspaces, we 
note 
\begin{proposition}\label{periodic}
Let $(S,\omega_S)$ be an $\aff \abx$ invariant subset.
Assume one of the coordinates $(z_h,z_v) \in S \subset X^2$ is 
periodic in $X$ (with respect to the action of $\aff \abx$).  
Then either $(S,\omega_S)$ is vertical, horizontal, 
or $(S,\omega_S)$ has nontrivial slope and 
both coordinates, $z_h$ and $z_v$ are periodic.
\end{proposition}
\begin{proof}
Assume $z_h \in X$ is periodic and $z_v$ 
is not, then 
\[\overline{\aff \abx \cdot(z_h,z_v)}=
\coprod_{z \in \aff \abx \cdot z_h} \{z\}\times X\]
is horizontal. If $z_v$ is periodic, but not $z_h$ 
we obtain the analogous statement with a 'vertical' subspace. 
Since a $2$-dimensional, invariant subspace $S$ with nontrivial slope 
cannot contain a vertical or horizontal subspace, 
$(z_h,z_v)\in S$ is periodic, or both coordinates $z_h$ 
and $z_v$ are non-periodic.  
\end{proof}
There is not much to say about 
the cone points of horizontal or vertical subspaces, 
since these spaces are isomorphic to a finite 
disjoint union of copies of $\abx$. 
\begin{proposition}\label{conepointorder}
Assume $(S,\omega_S)$ ($S \subset X^2$) is invariant, 
has nontrivial slope and $p=(p_h,p_v) \in S$ is periodic 
then 
\begin{equation}o_p(\omega_S)= 
\lcm(o_{p_h}(\omega),o_{p_v}(\omega))
\end{equation}
where $o_{p_h} = ord(p_h) +1$ (and $o_{p_v}$ is defined similarly).

\end{proposition}
\begin{proof}
We consider a full loop $\gamma \subset S$ around $p =(p_h,p_v)\in S$.  
Then $\gamma$ projects to two loops $\gamma_h \subset X_h$, 
$\gamma_v \subset X_v$ around $p_h$, $p_v$ respectively. 
The total angle along these two loops is a multiple of  
$o_{p_h}(\omega)$, $o_{p_v}(\omega)$ respectively. 
Thus the minimal total angle needed  
to return to the points $l_h \in \gamma_h $ and $l_v \in \gamma_v$ 
on the projected loops,   
while at the same time returning to  $(l_h,l_v) \in \gamma$, 
must be $\lcm(o_{p_h}(\omega),o_{p_v}(\omega))$. 
\end{proof}
\vspace*{1mm}


\section{Geometry of invariant subspaces.}

We give some examples of invariant subspaces in the arithmetic case.
\medskip\\
Given an invariant subspace $S \subset X_h \times X_v$ 
of slope $\alpha \in \R-\{0\}$, we assume in the following, 
that all images $pr_{h}(\CC_{\omega_S})\subset X_h$ and  
$pr_{v}(\CC_{\omega_S})\subset X_v$ of the set 
of cone-points $\CC_{\omega_S}$ of $(S,\omega_S)$ are {\em marked}. 
The leaves contained in a foliation of any 
subspace $S$ induce {\em affine maps} of leaves  
\begin{equation}
\begin{array}{ccccc}
\CF_{\theta}(X_h) & \stackrel{\quad pr_{h \ast}}{\longleftarrow}  & 
\CF_{\theta}(S)& 
\stackrel{pr_{v \ast}}{\longrightarrow}   & 
\CF_{\theta}(X_v)\\[2mm]
\CL_h            &\longmapsto      & \CL  &     \longmapsto &\CL_v
\end{array}
\end{equation}
which locally is a stretch by $\alpha$, because $S$ has slope $\alpha$. 
It turns out that the interesting sub-spaces have rational slopes.  
\begin{proposition}
Assume $S \in \CF_{\alpha}(X^2)$, with $\alpha=p/q \in \Q$, 
then the length of two compact, corresponding 
leaves $\CL_h \in \CF_{\theta}(X_h)$ 
and $\CL_v \in \CF_{\theta}(X_v)$ are related by: 
\begin{equation}\label{singrelation}
p\cdot |\CL_h|= q\cdot |\CL_v|, \quad 
\text{for singular leaves (saddle connections)}  
\end{equation} 
and there exist $a,b \in \N$
\begin{equation}\label{regrelation}
ap \cdot |\CL_v|= bq \cdot |\CL_v|, \ (a,b)=1
\quad \text{for compact regular leaves.}
\end{equation}
\end{proposition}
\begin{proof}
The proof uses the length induced by the $1$-forms on $S$ and $X$, 
i.e.
\begin{equation} |\CL| = \sqrt{1+\alpha^2 } \cdot |\CL_h|=
\frac{\sqrt{1+\alpha^2 }}{|\alpha|} \cdot |\CL_v| \quad \text{or, } \ 
\frac{|\CL_v|}{|\CL_h|}=|\alpha|.
\end{equation}
The condition $(a,b)=1$ comes from the fact that one 
completes one loop around $\CL \subset S$, starting let's say at 
$(x_h,x_v) \in \CL \subset S$,  if and only if one returns to 
$x_h \in \CL_h=\pr_{h \ast}\CL$ and  
$x_v \in \CL_v=\pr_{v \ast}\CL$ simultaneously while walking along $\CL$. 
\end{proof}
\vspace*{2mm}

\noindent{\bf Characterization of cone points in leaves of $\mathbf{\CF_{\alpha}(X^2)}$.}
Again we assume $\abx$ is a lattice surface. To describe $\slv{X,\omega}$ 
invariant subspaces of slope $\alpha \in S^1$ in $X^2$ is equivalent to finding 
the compact leaves of the foliation $\CF_{\alpha}(X^2)$. 
Here we study the neighborhood of a cone point 
$p=(p_h,p_v)\in X^2$, and the leaves of $\CF_{\alpha}(X^2)$ 
containing the cone point $p$.  
\smallskip\\
Two points $p_h,p_v \in X$ define a {\em cone point} 
$p=(p_h,p_v)$ of $X^2$, if at least one of the two points is 
a cone point of $\abx$. To study this cone point as a 
cone point of leaves in $\CF_{\alpha}(X^2)$ we take 
a small circle $\gamma_{r} \subset S$ of radius $r >0$ 
centered at $p$, $S$ a leaf in $\CF_{\alpha}(X^2)$ containing $p$, 
such that $p \in S$.  
The projections $\gamma_{r_h,\alpha} =\pi_h(\gamma_r) \subset X_h$ 
and $\gamma_{r_v,\alpha} =\pi_v(\gamma_r) \subset X_v$ of $\gamma_{r}$ 
are circles of radius $r_h=r/ \sqrt{1 + \alpha^2}$, 
$r_v=r\cdot|\alpha|/ \sqrt{1 + \alpha^2}$ respectively   
centered at $p_h$ and $p_v$. 
The product of the loops $\gamma_{r_h,\alpha}$ and $\gamma_{r_v, \alpha}$ 
is a torus 
\[\CT_{r, p, \alpha} := \gamma_{r_h,\alpha} 
\times \gamma_{r_v,\alpha} \subset X_h \times X_v. \]
We need to see how many leaves $ S \subset \CF_{\alpha}(X^2)$ 
terminate at $p$ and intersect the torus $\CT_{p,r,\alpha}$. 
A natural parametrization of $\CT_{p,r,\alpha}$, 
by means of the angle around $p_h$ and $p_v$ 
is $[0,o_{p_h}]\times [0,o_{p_v}]$. 
In these coordinates  with $*= h$ or $v$
\[\CF_{\theta}(X) \cap \gamma_{r_{\ast}} = 
\{\theta+i\ (o_{p_{\ast}}): i=1,...,o_{p_{\ast}} \}
.\] 

Consider the subspace 
$\CU_{\epsilon,p,\alpha} \subset X^2$ given by 
\[\CU_{\epsilon,p,\alpha}=\{p\} \cup 
\bigcup_{0< r <\epsilon} \CT_{p,r,\alpha}.\]
Geometrically it is a solid torus where the center loop is collapsed 
into one point, $p$.  It intersects all leaves of $\CF_{\alpha}(X^2)$ 
which contain $p$. We have an induced foliation 
\[  \CF_{\alpha}(\CU_{\epsilon, p,\alpha})
:=\CU_{\epsilon, p,\alpha} \cap 
\CF_{\alpha}(X^2).\]
The leaves of  $\CF_{\alpha}(\CU_{\epsilon, p,\alpha})$
containing $p$ are represented by
\begin{equation}\label{cyclic}
\begin{split}S^{\epsilon}_{i_h,i_v}:=\{(r_h\exp(2\pi i (\theta+i_h)/o_{p_h}), 
r_v \exp(2\pi i (\theta+i_v)/o_{p_v}):\\ 
\theta \in \R, 0 \leq r< \epsilon\} 
\subset \CF_{\alpha}(\CU_{\epsilon, p,\alpha}), 
\end{split}
\end{equation}
with $i_h=1,...,o_{p_h}$ and $i_v=1,...,o_{p_v}$. 
\medskip\\
To see how many of the $o_{p_h}\cdot o_{p_v}$ leaves 
$S^{\epsilon}_{i_h,i_v}$ are really different 
we derive from the representation (\ref{cyclic}) that 
\[S^{\epsilon}_{i'_h,i'_v}=S^{\epsilon}_{i_h,i_v}\]
whenever there is a $\theta \in \Z$, such that
\[i'_h \equiv \theta + i_h \mod o_{p_h}\ 
\text{ and }\ i'_v \equiv \theta + i_v \mod o_{p_v}.\]
Thus a necessary condition for the two local leaves to be 
equal is 
\begin{equation} i_p:\equiv i_h-i_v \equiv 
i'_h-i'_v \mod{\gcd(o_{p_h},o_{p_v})}
\end{equation} 
and given $i_p \mod \gcd(o_{p_h},o_{p_v})$,  
there are $\lcm(o_{p_h},o_{p_v})$ leaves $S^{\epsilon}_{i_h,i_v}$ 
such that $ i_h-i_v \equiv i_p \mod \gcd(o_{p_h},o_{p_v})$.
\begin{lemma}\label{cone-point-order}
Given two cone points $p_h,p_v \in X$. Then for every  
$\alpha \in \R-\{0\}$ the point $p=(p_h,p_v) \in X^2$ defines  
\begin{equation}
\frac{o_{p_h}\cdot o_{p_v}}{o_{p}}=\gcd(o_{p_h},o_{p_v})
\end{equation} 
cone points of order $o_p=\lcm(o_{p_h},o_{p_v})$ 
contained in leaves of $\CF_{\alpha}(X^2)$. 
In particular for small $\epsilon$ there are $\gcd(o_{p_h},o_{p_v})$
leaves $S^{\epsilon}_{i_h,i_v} \in \CF_{\alpha}(\CU_{\epsilon,p,\alpha})$.
\end{lemma}
\begin{proof}[{\bf Proof.}]
Recall that $o_p=\lcm(o_{p_h},o_{p_v})$ and note that 
the different cone points are defined as cone points 
of the local leaves $S^{\epsilon}_{i_h,i_v}$ containing them.   
Algebraically each cone point is characterized by the image of 
the map 
\begin{equation}\nonumber
\begin{array}{ccccc}
\Z/o_p\Z &\rightarrow &\Z/o_{p_h}\Z \oplus \Z/o_{p_v}\Z &\rightarrow& 
\Z/\gcd(o_{p_h},o_{p_v})\Z\\
i & \mapsto & (i_h,i_v):=(i (o_{p_h}),i (o_{p_v}))& 
\mapsto & i_h-i_v \ (\gcd(o_{p_h},o_{p_v})).
\end{array}
\end{equation}
Since for any given $i_p \in \Z/\gcd(o_{p_h},o_{p_v})\Z$ 
there is a leaf $S^{\epsilon}_{i_h,i_v}$ with 
$i_p \equiv i_h-i_v \mod \gcd(o_{p_h},o_{p_v})$ and 
\[o_h\cdot o_v=\lcm(o_{p_h},o_{p_v})\gcd(o_{p_h},o_{p_v})\] 
the statement follows.
\end{proof}
Alternatively one can characterize a cone point 
$p=(p_h,p_v) \in S \subset X^2$ 
by using the $o_h\cdot o_v$ images \
$\CL_{i,h}:=\pr_h(\CL_i)\subset \CF_h(X)$ and 
$\CL_{i,v}:=\pr_h(\CL_i)\subset \CF_h(X)$ of the horizontal leaves 
$\CL_i\subset \CF_h(S)$ emanating from $p$. For example 
the diagonal $D_+=\{(x,x)\in X^2: x \in X\}$ contains 
the cone point $p\in Z(\omega_{D_+})$, with singular 
leaves $\CL_{i}, i=1,...,o_p$ projecting to pairs of the shape  
$(\CL_{i,h}, \CL_{i,v})$.
\medskip\\
To give an example of an 
invariant surface $D_2 \subset \CF_1(X^2)$, $D_2 \neq D_+$ 
containing the cone point $p=(p_h,p_v) \in X^2$ ($p_h=p_v$ as cone point of $X$), 
we need to establish a correspondence of horizontal saddle 
connections starting at $p_h$ and $p_v$ in a cyclic way. 
To obtain a compact surface $D_2 \subset \CF_1(X^2)$, the easiest condition 
to impose on $X$ is all horizontal (and vertical) saddle connections have 
the same length. 
One arithmetic surface with this property is $L_3$, 
the $L$-shaped surface in genus $2$, tiled by three unit squares.  
\begin{figure}[h]
\epsfxsize=8truecm
\centerline{\epsfbox{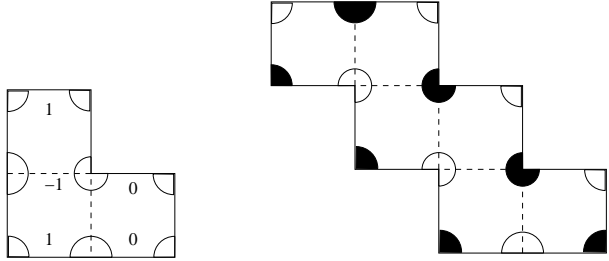}}
\caption{Slope $1$ invariant subspaces of $ L^2_3$ through $p$}
\refstepcounter{figure}\label{lshaped}
\end{figure}
\vspace*{1mm}\\
{\bf Note:} all horizontal and vertical relative periods of $D_2$ 
must have length $1$, by the cyclic identification scheme 
of their projections to $L_3=L_{3,h}=L_{3,v}$.  
The {\em staircase} shaped surface $D_2$ has an obvious automorphism $\phi$ 
of order $3$, which is the number of {\em steps} in the staircase. 
The map $\phi$ might be defined by moving each step of the staircase 
one step up. 
The quotient surface $D_2/<\phi>$ is a $2$-marked torus  
with Veech-group $\Gamma_2$.
\begin{theorem}
All $2$-dimensional $\slv{L_3,\omega_{L_3}}$-invariant subspaces 
of $S \subset L_3 \times L_3$ containing the cone point 
$p=(p_h,p_v)\in L^2_3$ are Veech surfaces $(S,\omega_S)$ 
with   
\[\slv{L_3,\omega_{L_3}} \cong \Gamma_2 \subseteq 
\slv{S,\omega_S} .\]
In particular $L_3 \times L_3$ contains infinitely many 
different arithmetic surfaces, whose Veech group contains the congruence group $\Gamma_2$.\vspace*{.5mm}  
\end{theorem}
\begin{proof}
The relative period lattice of $L_3$ defines a map 
\begin{equation}\label{torus-projection}
\pi: \slz \circ_d L^2_3 \longrightarrow \T^2\times\T^2 \cong \T^4.
\end{equation}
By construction $\pi$ is $\slz$ equivariant. Since the Veech-group 
of $L_3$ is $\Gamma_2$, $\slz \circ_d L^2_3$ has 
$3$ connected components and $L^2_3 \subset \slz \circ_d L^2_3$ 
as well as $p=(p_h,p_v)\in L^2_3$ is stabilized by $\Gamma_2$. 
By $\slz$ equivariance of $\pi$ all $2$-dimensional, 
$\slz$-invariant subspaces 
in $\T^4$ have $\Gamma_2$ invariant preimages in $L^2_3$. 
Theorem \ref{veech-groups} says that each connected 
component $(S,\omega_S)\subset L^2_3$ of an invariant 
subspace is stabilized by a group commensurable to $\Gamma_2$.
\medskip\\
By the discussion before Lemma \ref{cone-point-order}, there 
are exactly $3$ local leaves containing the point $p=(p_h,p_v)$. 
On each of these local leaves $p$ represents a cone point 
of order $3$. Take the horizontal foliation on $L_3$ and 
enumerate the $3$ horizontal leaves emanating from 
the cone point $p_h=p_v \in L_3$ by $(-1,0,1)$ 
(see enumeration in figure \ref{lshaped}), then the three cone 
points contained in the local leaves of slope $1$ in $L^2_3$ 
are characterized by the identification {\em types} 
$0 \leftrightarrow 0$, $0 \leftrightarrow 1$ and $0 \leftrightarrow -1$.  
\medskip\\
Certainly the leaf containing the cone point given 
by $0 \leftrightarrow 0$, i.e. the trivial identification, 
is stable under operation of $\Gamma_2$, because 
this cone point does not change type under the 
diagonal action of $\Gamma_2$. The other two cone 
points are getting exchanged under the involution $\phi$  
of $L_3$, because $\phi$ changes the order 
of outgoing leaves, i.e. $\phi (-1,0,1)=(1,0,-1)$. 
\medskip\\
Claim: the invariant subspace containing these two 
cone points is always connected. To see this note 
that any invariant subspace $S$ has rational slope, 
say $p/q$ $(p,q)=1$. Walking on a    
 horizontal leaf emanating from 
$(p_h,p_v) \in L^2_3$ on $S$, we will hit  
the cone point $(p_h,p_v) \in L^2_3$ again, since $S$ 
is a Veech surface. 
Choosing a continuation of this saddle connection gives a 
path $\gamma$ which is a chain of saddle connections 
and projects to 
a chain of $q$ (horizontal) saddle connections of 
length $1$ on $L_3 \cong  L_3 \times \{p_v\}$, respectively   
a chain of  $p$ saddle connections of the same length 
on  $L_3 \cong  \{p_h\} \times L_3$. This implies 
the middle point of $\gamma \subset S$ is either of the 
shape $(p_h,w)$, $(w,p_v)$, or of the shape $(w,w)$, 
$w$ a (regular) Weierstrass point. 
But any point of the above shape is fixed under the involution
$\phi$. This proves connectedness of $S$.     
As a consequence we have for the 
Veech-group of $(S,\omega_S)$: 
\[\Gamma_2 \subseteq \slv{S,\omega_S} .\] 
\end{proof}
\noindent {\bf Remark.} The leaf $\CL_2 \subset \CF_{2}(L^2_3)$ 
containing the standard neighborhood of the cone point $p \in L^2_3$ 
is tiled by $36$ squares, hence contains $12$ copies of $L$. 
It is easy to see that $\CL_{2}$ contains all $3$ cone points 
defined by $p=(p_h,p_v) \in X^2$. 
Also note, that two cone points are exchanged by 
the (only) involution $\phi \in \aff(L_3,\omega_{L_3})$. 
\medskip\\
{\bf Questions.} Describe the possible subspaces of $X^2$ 
for an arithmetic $\abx$ with one or more cone points. 
For given rational slope $\alpha$ and given cone point $p \in X$: 
how many leaves $\CL \in \CF_{\alpha}(X^2)$ contain $p$  
and how many squares are needed to tile $\CL$? What are the geometric 
and dynamic properties of $\CL$? How are these properties 
related to the properties of the marked surfaces 
parametrized by $\CL$. Ask the same question for 
a modular fiber  $\sF$ covering $X^2$, i.e. $\sF$ 
parametrizes branched covers of $X$. Can we 
say something about connectedness of $\sF$?   
%

\section{Invariant subspaces of $\T^4$}

If $\abx=(\R^2/\Z^2,dz)=:\T^2$ there are infinitely 
many $\slz$ invariant, complex subspaces in every  
dimension ranging from $0$ to $n$ in $\T^{2n}$. 
The description of these subspaces was (implicitly) given in \cite{s1}. 
From the geometric point of view all embedded unions of 
complex-linear tori  which define subgroups are $\slz$ invariant subspaces of $\T^{2n}$. 
We do not describe the higher dimensional case, i.e. subspaces of 
$X^n$ at this place, because it's a straightforward generalization 
of the $2$ complex dimensional case.   
\smallskip\\
\begin{defin}
We define a $\slz$-action on $\T^4$ ,by
\begin{equation}
\begin{array}{ccc} 
\slz \times\T^4 &\stackrel{\varodot}{\longrightarrow} 
& \T^4\cong \T^2\times\T^2\\
  (A,z)=(\left(\begin{smallmatrix} a & b \\ c & d 
\end{smallmatrix}\right),(z_h,z_v)) &\longmapsto& 
A\varodot(z_h,z_v):=(az_h+bz_v, cz_h+dz_v) 
\end{array}
\end{equation}
where $A=\left(\begin{smallmatrix} a & b \\ c & d 
\end{smallmatrix}\right)\in \slz$.
\end{defin}
Because the $\varodot$-action is complex-linear, it {\em commutes} with the 
diagonal-action of $\slz$ on $\T^4$. Indeed with $z_h:=x_{h}+iy_{h}$ 
and $z_v:=x_{v}+iy_{v}$ one easily verifies 
\begin{equation}
\left(\begin{smallmatrix} e & f \\ g & h 
\end{smallmatrix}\right)\varodot \left[ \left(\begin{smallmatrix} a\ b \\ c\ d 
\end{smallmatrix}\right)\circ_d \left(\begin{smallmatrix} x_{h}+iy_{h} \\x_{v}+iy_{v}
 \end{smallmatrix}\right)\right]=\left(\begin{smallmatrix} a\ b \\ c\ d 
\end{smallmatrix}\right) \circ_d \left[ \left(\begin{smallmatrix} e & f \\ g & h 
\end{smallmatrix}\right)\varodot 
\left(\begin{smallmatrix} x_{h}+iy_{h} \\x_{v}+iy_{v}
 \end{smallmatrix}\right) \right]
\end{equation} 

\noindent We call an $\slz$ invariant subspace $(S, \omega_S) \subset \T^4$ 
{\em simple}, if there is a $z \in S$ with 
$S=\overline{\slz \cdot z}$, i.e. if $S$ is the closure of the $\slz$ 
orbit of a single point. 
\smallskip\\
Ratner's Theorem implies that the that $\slz$-invariant subspaces 
of the Lie-group $\T^4$ are Lie-subgroups. Thus the classification 
of $\slz$-invariant subspaces is well known and easy to achieve. 
The purpose of this section is to characterize invariant subspaces  
using the slope defined earlier. 
To state the theorem it is useful to write
\[\sO_n := \left\{\left[\frac{a}{n},\frac{b}{n}\right]: a,b,n \in 
\Z^2, (a,b,n)=1\right\}
=\slz \cdot \left[\frac{1}{n},0\right] \subset \T^2.\] 
\begin{theorem}\label{torus-subs}
The following statements are equivalent
\begin{itemize}
\item[1.] $S \subset \T^4$ is a $2$-dimensional, simple  
$\slz$ invariant subspace
\item[2.] Given $S \subset \T^4$, there is an $A \in \slz$ and an $n \in \N$, 
such that $S_h:=A\varodot S$ 
is horizontal and $S_h \cap \T^2_v =\sO_n$.
\item[3.] There is a nontrivial, $\slz$-equivariant  
homomorphism of complex Lie groups $\pi: \T^4\cong \T^2 \times \T^2 \rightarrow \T^2$,  
which does not factor over an isogeny and an $n \in \N$, 
such that $S=\pi^{-1}(\sO_n)$. 
\item[4.] $S$ consists of solutions of $az_h+bz_v \equiv 0 \mod \Z^2 $ 
with $a,b \in \Z$ and $(a,b)=n$ which do not solve any equation 
with coefficients $(a,b)=m$ for $m|n$.
\end{itemize}
\end{theorem}
\begin{proof}
From Theorem \ref{real} we know that all invariant 
subspaces are defined by a real, linear equation. 
Hence all invariant subspaces 
are contained in the foliations $\CF_{\alpha}(\T^2 \times \T^2)$,  
$\alpha \in \R \cup \{\infty\}$. 
It is well known that leaves of $\CF_{\alpha}(\T^2 \times \T^2)$ 
with irrational slope are dense in $\T^4$. 
Thus any closed, invariant $S \subset \CF_{\alpha}(\T^2 \times \T^2)$ 
is contained in $\CF_{\alpha}(\T^2 \times \T^2)$, with $\alpha \in \Q \cup \{\infty\}$. 
But for rational $\alpha$ there is an $A \in \slz$, 
such that $A\varodot \CF_{\alpha}(\T^2 \times \T^2)=\CF_{h}(\T^2 \times \T^2)$.  
Since the real-action commutes with the diagonal action, 
we know $S_h := A\varodot S \in \CF_{h}(\T^2 \times \T^2)$ 
is $\slz$ invariant and for the same reason $S_h$ is simple, if $S$ is. 
The  $\slz$-orbit classification on $\T^2$ 
then implies there is an $n \in \N$
such that $S_h \cap \T^2_h=\sO_n$. Obviously $(2)$ implies $(1)$,  
showing the equivalence of statement $(1)$ and $(2)$. 
Note, that each component of $S_h$ is a torus $\T^2$, which is indeed 
the orbit closure of a single (irrational) point under the action of $\slz$. 
\medskip\\
$(2) \Leftrightarrow (3)$ Take $A \in \slz$ making $S$ horizontal 
and look at the map 
\begin{equation}
\psi: \begin{array}{cclc}
 \T^4 &\longrightarrow& \T^2\times\T^2 &\longrightarrow \T^2\\
z &\longmapsto& (z_h,z_v):=A\varodot  z &\longmapsto z_v
\end{array}.
\end{equation}
Since both maps in the composition are linear and commute with the (diagonal) 
action of $\slz$, $\pi$ is a $\slz$-equivariant homomorphism 
of the complex Lie groups $\T^2 \times \T^2$ and $\T^2$. 
Now $\pi$ cannot factor over an isogeny 
\[ \psi: \T^2 \rightarrow \T^2, \quad [z] \mapsto [az]  \quad \text{ where } a \in \Z,  \] 
because $\det A =1$.
The image of $\pi$ is by construction $\sO_n$ for some $n \in \N$. 
On the other hand every complex Lie group homomorphism
 $\pi:\T^2 \times \T^2  \rightarrow \T^2$ is given by 
 \[(z_h,z_v)\mapsto cz_h+dz_v \quad \text{with } c,d \in \Z.\]
If $\pi$ does not factor over an isogeny  we must have $(c,d)=1$ 
and we can extend $\pi$ to a linear map 
\[ A:\T^2 \times \T^2 \rightarrow \T^2_h \times \T^2_v, \quad A \in \slz\]
where $A=\left[\begin{smallmatrix} a\ b \\ c\ d\end{smallmatrix}\right]$. 
Now the image of $S=\pi^{-1}(\sO_n)$ under $A$ is horizontal 
and simple because $\sO_n \subset \T^2$ is.     
\medskip\\
For $(2) \Rightarrow (4)$ we note that any $A =
\left[\begin{smallmatrix} a\ b \\ c\ d\end{smallmatrix}\right]\in \slz$   
making $S$ horizontal provides us with $\sO_n$ equations characterizing $S$
\begin{equation}\label{oneq} 
cz_h+dz_v  \equiv \left[\frac{p}{n},\frac{q}{n}\right] \mod \Z^2, \quad \text{where}  
\left[\frac{p}{n},\frac{q}{n}\right]  \in \sO_n.  
\end{equation}
These can be expressed by the equation
\[ ncz_h+ndz_v  \equiv 0 \mod \Z^2.\]
Now the statement follows because $(c,d)=1$. 
For the reverse implication we easily see that given an equation 
$ncz_h+ndz_v  \equiv 0 \mod \Z^2$  with $(c,d)=1$,  
the set of solutions which solves no equation of the shape 
$mcz_h+mdz_v  \equiv 0 \mod \Z^2 $ for any $m|n$ 
solves precisely the $|\sO_n|$ equations (\ref{oneq}). 
Using $c$ and $d$  as second row in a matrix $A \in \slz$ in 
turn defines a map which makes the set of solutions 
of the given equation horizontal and has the desired intersection 
property.
\end{proof}

\begin{proposition}\label{t4inv}
Let  $Y$ be a subset of  $\T^4$ which is closed, invariant under the diagonal
action of  $\slz$ and contains infinitely many invariant two dimensional torii then $Y$ is equal to $\T^4$.
\end{proposition}

We first prove an arithmetical lemma:

\begin{lemma}\label{equispace}
Given any $\varepsilon >0$, the orbit 
\[
\sO_n = \slz \cdot \left[\frac{1}{n},0\right] \subset  
\T^2
\]
is $\varepsilon$ dense in $\T^2$ for $n$ large enough.
\end{lemma}

\begin{proofof}{Lemma \ref{equispace}}
The set $\sO_n$ contains the set
$$\left\{\left[\frac{a}{n},\frac{b}{n}\right]: a,b,n \in 
\Z,\  \gcd(a,n)=1 \mbox{ and } \gcd(b,n) = 1 \right\}.$$
Therefore, it is enough to prove that the set 
$\left\{\frac{a}{n}: a,n \in 
\Z, (a,n)=1\right\}$
is $\e$ dense in $\T^1$ for $n$ large enough.

The Jacobsthal's function $J(n)$ is the largest gap between consecutive integers relatively prime to $n$ and less than $n$. 
We have to prove that $J(n)/n$ tends to 0 as $n$ tends to infinity to complete the proof of lemma \ref{equispace}.
In fact, a much more precise result does exist.
By Iwaniec's result (see \cite{iv}) there is a constant $K$ such that 
$J(n) \leq K (\ln n)^2$.
This ends the proof of lemma \ref{equispace}.
\end{proofof}

\begin{proofof}{Proposition \ref{t4inv}}
 We use the classification of the orbit closures obtained
in part 4) of Theorem \ref{torus-subs}:  for the $r$th torus $S_r$
we denote the coefficients of the equation of the torus
by $a_r,b_r,n_r$.  If the coefficients are uniformly bounded (in $r$)
then there are only a finite number of tori, thus at least one
of the coefficients is not bounded.  Given $\e > 0$,
if $n_r$ is sufficiently large then 
the orbit of the torus $S_r$ (by Theorem \ref{torus-subs} part 2)  
is $\e$ dense since it consists of 
$\varphi(n)\psi(n)$ parallel
tori with transverse spacing at most $\e$  by Lemma \ref{equispace}.

On the other hand, if $n_r$ is bounded then the slope sequence $a_r/b_r$ has a
limit point (maybe 0 or infinity).  If this limit ``slope'' is
irrational then clearly $Y = \T^4$.  In the rational case we also have $Y = \T^4$, since the 
area of the $S_r$ approach
infinity with transverse spacing approaching $0$.

\end{proofof}

\section{Proof of Theorem \ref{finite} and Theorem 
\ref{classification}}

\noindent With the previous observations we are ready to prove 
Theorem \ref{finite}:
\begin{proof}[{\bf Proof of Theorem \ref{finite} for $\mathbf{X^2}$}]
If $X$ is arithmetic then the invariant subspaces are lifts of the
invariant subspaces of $\T^4$.  Since there are infinitely many
invariant subspace of $\T^4$ (see Theorem \ref{torus-subs}) there are infinitely
many for $X$.

Let $S$ be a 2-dimensional, $\aff \abx$ 
invariant submanifold of $X^2$ with nontrivial slope. 
Then the two projections 
\begin{equation} \begin{CD}
X_h @<{\pr_h}<< S @>{\pr_v}>>  X_v 
\end{CD} 
\end{equation}
are surjective and 
every point $p=(p_h, p_v) \in S$ with one periodic coordinate 
is already periodic, by Proposition (\ref{periodic}). 
If $\abx$ is not arithmetic, it contains only finitely 
many periodic points (see \cite{ghs}).  
Thus in the non-arithmetic case, we can look at 
the surface $(X,p_1,...,p_n)$, 
marked in all the periodic points 
(some of them might be cone points). 
Note that by assumption the set of periodic points is 
not empty, because $g(X)\geq 2$ for non-arithmetic $\abx$. 
\medskip\\
Any completely periodic 
foliation $\CF_{\theta}(X)$ on $(X,p_1,...,p_n)$ contains 
a nonempty {\it spine}, i.e.~the weighted graph of saddle connections 
contained in $\CF_{\theta}(X)$, the weights being the lengths of the saddle connections.  
The spine is nonempty by the existence of cone points. 
With respect to the two projections above 
any invariant subspace $S$ defines a {\it similarity} of spines 
from $X$ to itself, i.e. a bijective map of 
the spines of $\CF_{\theta}(X_h)$ and $\CF_{\theta}(X_v)$
stretching the length of each saddle connection by 
a factor $|\alpha|$, $\alpha \in \R^{\ast}$ the slope 
of $S$.
\medskip\\
Without restrictions we can assume    
$\alpha \in \R^{\ast}$ satisfies $\vert \alpha \vert \geq 1$, otherwise we change 
the roles of $pr_h$ and $pr_v$. 
Now take the longest saddle connection in the horizontal 
spine of $X$. If $|\alpha|>1$ this saddle 
connection would be mapped to a saddle connection stretched 
by a factor $|\alpha|$, contained in 
the horizontal spine of $X$. Contradiction, thus 
we must have $\alpha=\pm 1$.
\medskip\\
For slope $\pm 1$ subspaces $S \subset X^2$ we first note that 
$S$ intersects the fiber $\{p_h\} \times X_v \subset X_h \times X_v$, 
where $p_h \in X$ is any cone point. To see this take a point 
$(z_h,z_v)\in S \subset X_h \times X_v$ and connect $z_h$ with $p_h \in X$ 
using a path in $X = X_h$, this path lifts to $S$. In other 
words, the two coordinate projections of $S \rightarrow X_i$, $i=h,v$ 
are surjective. By Proposition \ref{periodic} the 
intersection $S \cap \{p_h\} \times X_v$ 
consists of periodic points, 
since $p_h \in X$ is periodic. Thus the 
finiteness of the number of subspaces follows by the finiteness of 
periodic points and finiteness of orders of cone-points on $X$ by lemma 
\ref{cone-point-order}. 
Note that any neighborhood $U \subset S \subset X^2$ ($S$ invariant) 
of a periodic point determines $S$ by ergodicity of the $\slv{X,\omega}$ 
action on $S$.    
\end{proof}
\begin{proof}[{\bf Proof of Theorem 
\ref{real} and \ref{finite} for $\mathbf{X^n, \ n \geq 3}$}]
We prove the claim by induction over the complex 
dimension $n$ of $X^n$. We already have the result for 
$n=1$ (finitely many periodic points) and for $n=2$. 
The claim is clear, if $X$ is arithmetic, thus we assume 
$X$ is not arithmetic. In this case we can also assume 
all periodic points of $X$ are fix points. This does 
not change the invariant subspaces $S \subset X^n$, 
since the group $\Gamma \subset \slv{X,\omega}$ 
stabilizing all periodic points of $X$ 
is a finite index subgroup of $\slv{X,\omega}$. 
Without restrictions we assume $\slv{X,\omega}$ stabilizes 
all periodic points. 
Denote by 
\[pr_i: X^n \rightarrow X^{n-1}_i:= X_1 \times ... 
\times X_{i-1} \times X_{i+1}\times ...\times X_{n}\]
the projections to the $n$ {\em hyperplanes} $X^{n-1}_i$. 
Each periodic($=$ fix) point $p\in X$ defines an embedding 
\[e_{i,p}: X^{n-1}_i \rightarrow  X_1 \times ... 
\times X_{i-1} \times \{p\}\times X_{i+1}\times ...\times X_{n} 
\subset X^n\]
of each hyperplane $X^{n-1}_i$ with the property 
\[\pr_i \circ e_{i,p} =\id_{X^{n-1}_i}, \ \text{for all } 
i=1,..,n.\] 
For any {\em invariant} hypersurface $S \subset X^n$ and any 
$\slv{X,\omega}$-fix point $p \in X$ we have  
\begin{equation} 
\dim_{\C}(e_{i,p}(X^{n-1}_i)\cap S)=\left\{ 
\begin{array}{l}
n-1 \\
n-2 \\
0
\end{array}
\right. .
\end{equation} 
In the first case $ \slv{X,\omega}\cdot e_{i,p}(X^{n-1}_i) \cong S$,  
while in the last case $e_{i,q}(X^{n-1}_i) \cong S$ 
for an $ \slv{X,\omega}$-fix point $q \neq p$. 
In the remaining case $S$ is characterized by $n$ intersections,  
which are invariant subspaces of lower dimension. In fact,  
choosing an $\slv{X,\omega}$-fix point $p \in X$ and considering  
the $n$ hyperplanes $X^{n-1}_i:=e_{i,p}(X^{n-1})$ shows that 
$S$ is defined by the $n$ intersections $S\cap X^{n-1}_i$.
\medskip\\ 
If $\dim_{\C}(S)\leq n-2$ we can describe $S$ as an intersection 
of invariant surfaces of dimension less than $n-1$: 
\begin{equation}
S=\bigcap_{i}\pr^{-1}_i\circ \pr_i(S).
\end{equation}
The statement follows now by induction over the dimension $n$.
\end{proof}
\begin{proposition}\label{diagonal}
Given a reduced lattice surface $\abx$, then 
the diagonal $(D_{+},\omega_+)$ and the off-diagonal (if it exists) 
$(D_{-},\omega_-)$, $ (D_{\pm} \subset X^2)$ are translation surfaces 
isomorphic to $\abx$. 
\end{proposition}
\begin{proof}
For the diagonal 
$\omega_{\pm}=\sqrt{2}\cdot \omega$ is induced by the bijective translation
\begin{equation}
\begin{array}{ccc}
D_+ & \longrightarrow & X \\
(x,x) & \longmapsto & x
\end{array}
\end{equation}  
while for the offdiagonal $w_- = \sqrt{2} \omega$ is induced by
\begin{equation}
\begin{array}{ccc}
D_- & \longrightarrow & D_+ \\
(x,\phi(x)) & \longmapsto & (x,x)
\end{array}
\end{equation} 
where $\phi \in \aff \abx$ is the unique 
map, such that $\D\phi = -\id$. 
\end{proof}
\begin{proof}[{\bf Proof of Theorem \ref{veech-groups}}]
All invariant subspaces $S \subset X^n$ of complex dimension $1$ are real-linear and $\aff \abx$ acts 
as a group of real-affine homeomorphisms with respect to the given 
complex-linear structure. We obtain the claim after noticing 
that each connected component of $S$ has a lattice stabilizer 
in $\aff \abx$, in particular this 
stabilizer is $\slv{X,\omega}$ if $S$ is connected. 
Note that as a Veech surface a priori $S$ 
could have an even bigger Veech- and affine group.  
\medskip\\ 
We explicitly define the differential $\omega_S$ on $S$ 
using local parametrizations of $S$ which are induced by 
natural charts of $\abx$. 
\medskip\\
By real linearity of $S$ 
there are vectors $(a_1,...,a_n)\in \R^n$,  
$(t_1,...,t_n)\in \C^n$
and $n$ local translation maps $S \rightarrow X$, 
such that   
\begin{equation}\label{transmap}
\phi: 
\begin{cases}
S\hookrightarrow S\subset X^n &\\
z \longmapsto (a_1z+t_1,a_2z+t_2,...,a_nz+t_n)&
\end{cases}
\end{equation}
defines a local homeomorphism. Now the  
$1$-form on $S \subset X^n$ defined by 
\[ \frac{1}{\|a\|^2}\sum^n_{i=1} a_i\omega_i|_{TS}, \ \text{ with }\ \|a\|^2:= \sum^n_i a^2_i\]
implies that locally 
\begin{equation}
\phi^{\ast}\omega_S=\phi^{\ast}\left(\frac{1}{\|a\|^2}\sum^n_{i=1} a_i\omega_i|_{\smash[b]{TS}}
\right)=\frac{1}{\|a\|^2}\sum^n_{i=1} a^2_i dz=dz
\end{equation}
showing that  (\ref{transmap}) defines a translation map with respect to 
natural charts.
The diagonal action of $\aff \abx$ on $X^n$ induces an 
action of $\slv{x,\omega}$ on $\omega_S$, given by
\[ A.\omega_S:= \frac{1}{\|a\|^2}\sum^n_{i=1} a_i A\cdot \omega_i|_{\smash[b]{TS}} 
=A \cdot \frac{1}{\|a\|^2}\sum^n_{i=1} a_i\omega_i|_{\smash[b]{TS}}= A \cdot \omega_S\]
The last identity is caused by real linearity of $A \in \slv{x,\omega}$. 
But this is the claim, i.e. the diagonal 
action induces the standard action of $\aff \abx$ on $(S,\omega_S)$ 
given by 
\begin{equation} A \cdot \omega_S := \left( 1, i \right)
\left(\begin{matrix}a & b \\ c & d \end{matrix}\right) 
\left(\begin{matrix} \re (\omega_S) \\ \im (\omega_S) \end{matrix}\right),  
\quad \text{ if } \quad A=\left(\begin{matrix}a & b \\ c & d \end{matrix}\right) 
\end{equation}
\end{proof}
\begin{lemma}\label{arithmetic}
Take a lattice surface $\abx$. Assume there is 
a periodic direction $\theta$, such that   
all absolute periods represented by leaves $\CL \in  \CF_{\theta}(X)$ 
fulfill $|\CL| \in \alpha \cdot \Q$, $\alpha \in \R_+$. 
Then $\abx$ is arithmetic. 
\end{lemma}
\begin{proof}
Without restrictions we can assume $\alpha = 1$, the periodic 
direction $\theta$ is horizontal and 
the vertical foliation of $\abx$ is periodic too. 
Take the horizontal cylinder decomposition of $\abx$, the  
$i$-th cylinder having height $h_i$ and width $w_i$. 
Since $\abx$ is 
a Veech surface we have $w_i/h_i \cdot h_j/w_j \in \Q$ 
and thus $h_i/h_j \in \Q$ for all pairs of cylinder(-heights), 
because all $w_i$ are rational.  
This implies 
\[h_1,...,h_n \in \alpha_h \cdot \Q \] 
and again we can assume $\alpha_h=1$. Since all the width and 
the heights of cylinders, including the absolute horizontal 
periods, which might represent cylinders of height $0$, 
are rational, they define a lattice in $\R^2$. 
This in turn defines a translation map 
$\pi: \abx \rightarrow (\R^2/\per (\omega),dz)$, 
proving the Lemma. 
Note that by assumption 
$\per(\omega) \subset \Z^2$. 

To represent $\pi$ (see also \cite{zo2}), 
choose a cone point $p \in X$ and any other point $z \in X$. 
Then any path connecting $p$ with $z$ is isotopic (in $X-Z(\omega)$)  
to a path consisting entirely of horizontal and vertical line  
segments. Now take the holonomy of the straightened path modulo 
$\per(\omega)$ and note that any other path from $p$ to $z$ 
gives the same vector modulo $\per(\omega)$. 
\end{proof}

\begin{proof}[{\bf Proof of Theorem \ref{classification}}]
First note, that for primitive and by definition for 
reduced lattice surfaces $\abx$ taking derivatives induces an 
isomorphism $\slv{X,\omega} \cong \aff \abx$. Thus we can assume 
$\slv{X,\omega}$ acts on $X^2$.
\smallskip\\
By Theorem \ref{torus-subs} the statements hold 
for arithmetic surfaces. Thus we might assume $\abx$ is not arithmetic 
and all the (finitely many) periodic points of $\abx$ are marked.  
\smallskip\\
Assume $S \subset X^2$ is a $\aff \abx$ invariant subspace 
of slope $\pm 1$. Because we mark all periodic points of $X$, 
$S$ induces a length preserving bijection of saddle connections  
for each completely periodic direction. 
Now take a cone point $p=(p_h,p_v)\in S \subset X_h \times X_v$, 
such that (without restrictions) $p_h$ is 
the (one and only) cone point of $X=X_h$. 
We have two cases:  
\smallskip\\
$p_v$ is not a cone point, but periodic. Then 
for any given periodic direction $\theta$, we consider 
a (chain of) saddle connection(s) $\CL_{p_v}\in \CF_{\theta}(X)$ connecting 
$p_v$ with itself, containing $p_v$ only once. 
Again there are two cases:  \\
$\CL_{p_v}\in \CF_{\theta}(X)$ is a regular leaf. 
Then $S$ sets up a correspondence of $\CL_{\theta}(p_v)$ 
and all $o_{p_h}$ leaves $\CL_{p_h, i}\in \CF_{\theta}(X)$ 
 $i=1,...,o_{p_h}$ starting and ending at $p_h$ must have 
length 
\begin{equation}\nonumber
|\CL_{p_v}|=|\CL_{p_h, i}| \quad \text{ for all } i=1,...,o_{p_h}
\end{equation} 
The other case is, if $p_h \in \CL_{p_v}$ $p_h=p$ {\em the} cone point 
of $\abx$. Again we require that $\CL_{p_v}$ meets $p_h$ just one 
time. Then $\CL_{p_v}=\CL_{p_h, i}$ for one $i$, say $i=1$ and 
via the correspondence defined by $S$ we find again  
\begin{equation}\nonumber
|\CL_{p_h, 1}|=|\CL_{p_h, i}| \quad \text{ for all } i=1,...,o_{p_h}.
\end{equation}
This finishes the first case. 
\medskip\\
In the second case  $(p_h,p_v)=(p,p)$ is the cone point of $\abx$ 
in both coordinates. Again we denote all $o_{p_h}$ singular leaves  
connecting $p=p_h$ with itself in a completely periodic direction 
$\theta$ by $\CL_{p_h, i}$, $i=1,...,o_p$. We know that locally 
there are $o_{p_h}$ $2$-dimensional leaves in $\CF_{\pm 1}(X^2)$ 
containing $(p,p)$, labeled by which leave $\CL_{p_h, i}$ we identify 
with $\CL_{p_h, 1}$. If $o_{p_h}$ is prime every nontrivial 
combination will again give 
\begin{equation}\nonumber
|\CL_{p_h, 1}|=|\CL_{p_h, i}| \quad \text{ for all } i=1,...,o_{p_h}
\end{equation}
by cyclicity of the identifications. 
The trivial identification, i.e. $\CL_{p_h, 1}$ with $\CL_{p_v, 1}$, 
is contained in $D_{\pm}$. 
Both cases imply with Lemma \ref{arithmetic}: 
either $\abx$ is arithmetic in contradiction to 
our assumption, or $S=D_{\pm}$.    
\end{proof}
Note: we prove more, than stated in Theorem \ref{classification}, 
in fact we can say
\begin{corollary}
Let $\abx$ be a reduced Veech surface with 
exactly 
\underline{one} cone-point $p \in X$, 
then all $2$-dimensional 
$\slv{X,\omega}$-invariant subspaces $S \subset X^2$ 
with non-trivial slope contain the cone point $(p,p) \in X^2$. 
\end{corollary}
\noindent The Corollary implies that $X^2$, $\abx$ a reduced, non-arithmetic 
lattice surface with exactly one cone point, 
cannot contain more than $o=o_p$ connected, invariant 
surfaces of slope $\pm 1$. It follows from Proposition \ref{periodic} 
and Lemma \ref{cone-point-order} 
that for non-arithmetic Veech surfaces $\abx$ the number 
of compact leaves (with slope $\pm 1$) is bounded by 
\[ \begin{split}\# \{S\subset X^2: \dim_{\C}(S)=1, \text{ slope}(S)=\pm 1\}\leq \\ 
\leq \frac{\displaystyle \sum_{p,q \in X \text{ periodic}}o_p\cdot o_q}
{\displaystyle \sum_{p \in X \text{ periodic}}o_p} 
=\sum_{p \in X \text{ periodic}}o_p, 
\end{split}
\]
which is the {\em total order} of periodic points on $\abx$. 
To make this bound effective one can use  Möller's result \cite{mm} 
for the number of periodic points on non-arithmetic Veech surfaces 
of genus $2$ with one cone point. 
We obtain the upper bound $8$ for the number of (off-) diagonals 
for genus $2$ surfaces with one cone point.
\medskip\\   
\noindent Another immediate consequence is 
\begin{corollary}\label{knownexample}
All primitive Veech surfaces $\abx$ of genus $g$, with 
exactly one cone point, such that $2g-1$ is prime have only $D_{+}$ 
and $D_-$ as $2$-dimensional, $\slv{X,\omega}$-invariant 
subspaces of $X^2$. $D_-$ exists if and only if $\slv{X,\omega}$ 
admits the element $-\id$. 
In particular for all primitive Veech surfaces 
$\abx$ with one cone point 
in genus $2,3$ and $4$ the only invariant 
surfaces are $D_{-}$ and $D_{+}$.
\end{corollary}
\begin{proof}
We have $o_p=2g-1$, if there is only one cone point $p$ 
on $\abx$ and  $o_p=3,5,7$ for $g=2,3,4$.
\end{proof}
Note that all known primitive Veech surfaces 
in genus $2,3$ and $4$ possess an involution, 
i.e. $-\id \in \slv{X,\omega}$.


\section{The illumination problem --- prelattice surfaces}

A discrete group $\Gamma \subset \slr$ called {\em prelattice}, 
if it contains noncommuting parabolic elements. 
A translation surface $\abx$ is a 
{\it prelattice surface}  if $\slv{X,\omega}$ is a prelattice group.
By acting by $\slr$ we see that 
if a surface is prelattice then we can always assume that the
two parabolic elements are orthogonal, one horizontal one 
vertical.\footnote{We make this assumption for simplicity of 
notation, it affects
none of our results.}
Thus we can decompose the surface into vertical cylinders and
also into  horizontal cylinders.  Taking the intersection of the
cylinders leads to a decomposition  into rectangles with vertical 
and horizontal sides.
For $p \in X$ let $\CC_v(p)$ be the vertical cylinder,
$\CC_h(p)$ the horizontal cylinder and $\Rec(p)$ the rectangle
containing $p$. 

Let 
$$T_h := \left ( \begin{matrix}
1 \ a \\
0 \ 1
\end{matrix} \right ) \quad \hbox{ and } \quad
T_v := \left ( \begin{matrix}
1 \ 0 \\
b \ 1
\end{matrix} \right )
$$
be the horizontal and vertical parabolic elements in $\slv{X,\omega}$.
Here $a$ and $b$ depend on the prelattice surface $\abx$.  
Suppose $\CC$ is a vertical cylinder with height $h$ and width $w$
and with origin the bottom left corner.
Then the action of $T_v$ on $\CC$ is given by a power of the
Dehn twist:
\begin{equation}\label{vangle}
T_v\left (  \begin{matrix} x\\y \end{matrix} \right ) =  
\left (  \begin{matrix} x\\y + bx \ ( h ) \end{matrix} \right ).
\end{equation}
Here $( h )$ means modulo $h$.
Similarly for a horizontal cylinder $\CC$ the action of $T_h$ is given by
\begin{equation}\label{hangle}
T_h\left (  \begin{matrix} x\\y \end{matrix} \right ) =  
\left (  \begin{matrix} x + a y\ ( h ) \\y  \end{matrix} \right ).
\end{equation}
Warning, in the standard notation which we follow here
the height of a horizontal cylinder is noted $w$ and the width is noted
$h$.
\medskip\\
If $\slv{X,\omega}$ is a lattice, 
then $\abx$ is a {\it Veech surface}.
Veech has shown that the direction of any saddle connection on a Veech surface
is a parabolic direction \cite{v1}. 
\begin{proposition}\label{inva}
The set of  $\{(p,q): p \hbox{ illumines } q\}$ is invariant under
the diagonal action of $\aff \abx$, i.e.~if $p$ illumines
$q$ then $\phi(p)$ illumines $\phi (q)$ 
for all $\phi \in \aff \abx$.
\end{proposition}
This follows from the fact that $\phi$ takes a geodesic to a geodesic.
\medskip\\
A key observation is the following simple fact.
\begin{proposition}\label{prop2}
Fix a vertical (or horizontal) cylinder $\CC \subset X$. 
Each point $p \in \CC$ illumines every $q \in \Rec(p)$.
\end{proposition}
This follows since a rectangle is convex.
Theorem \ref{thm1} follows immediately
from Propositions \ref{inva} and \ref{prop2} combined with
the following more general theorem.
\begin{theorem}\label{thm2}
Suppose that $\abx$ is a prelattice surface and 
fix $p \in X$.  Then except for an at most countable exceptional set of
$q \in X$ there exists a $\phi := \phi(p,q) \in \aff \abx$
such that $\phi(p)$ and $\phi(q)$ are both contained in
the rectangle $\Rec(p)$.
\end{theorem}

Part of the proof of Theorem \ref{thm2} is based on the arithmetic
analysis of rotations of the torus.  We consider the two torus
$\T^2 = \S^1 \times \S^1$ and for $(c,d) \in \T^2$ let
$R^n_{\varphi,\vartheta}(c,d)$ be the rotation by $(\varphi,\vartheta)$.
We need the following quantitative version of Kronecker's theorem.

\begin{lemma}\label{Kronecker}
Fix $(c,e) \in \S^1 \times \S^1 \cong \T^2$, $\varphi \in \S^1$ and $\e >0$.  
Then for all, but finitely many $\vartheta \in \S^1$ and for 
all $d \in \S^1$ there exists an $n$, such that 
$R^n_{\varphi,\vartheta}(c,d) \in (c - \e,c +\e) \times (e-\e,e+\e)$.
\end{lemma}

{\em Remark}: 
the cardinality of the exceptional set depends on $\varphi$ and $\e$.

\begin{proof}
Consider a $\vartheta \in \S^1$. There are thee cases according to the
dimension $\Delta$ of the vector space over $\Q$ 
generated by $(1,\varphi,\vartheta)$.

Suppose this dimension is $3$.
The classical Kronecker theorem asserts that 
for all $d$, the orbit of $(c,d)$ by $R^n_{\varphi,\vartheta}$ is dense in
$\T^2$ (\cite{KaHa}, page 29) and in this case the exception set of
$\vartheta$ is empty.
\medskip\\
Next suppose the dimension satisfies $\Delta = 1$.
Thus $\varphi$ and $\vartheta$ are rational numbers, we note them in
reduced form as 
$\varphi = p_1/q_1$ and $\vartheta = p_2/q_2$.
We emphasize that the number $\varphi$ is fixed.  For any $n$ the point 
$R^{nq_1}_{\varphi,\vartheta}(c,d)$ has first coordinate $c$.
Thus we think of $R^{nq_1}_{\varphi,\vartheta}(c,d)$ as a rotation of
the circle $\{c\}\times \S^1$.
Thus if $\theta$ is rational with denominator $q_2$ 
greater than $1/\e$, 
then orbit $R^{nq_1}_{\varphi,\vartheta}$ of any point $(c,d)$ is 
$\e$-dense in the circle $\{c\} \times \S^1$ and the lemma holds with
the exceptional set of $\vartheta$ (a subset of) 
those with denominator less than $1/\e$. 

Finally consider the case when $\Delta = 2$. 
Thus there exists $A,B,C \in \Z$  with $\gcd(A,B,C) =1$
such that $A\varphi + B\vartheta = C$. Let $D := \gcd(A,B)$,
$A' := A/D$ and $B' :=B/D$.
Then the above equation becomes $A' \varphi + B'\vartheta = C/D$.

By homogeneity it suffices to prove the lemma for $c=d=0$.
The orbit of $(0,0)$ is contained in the torus-loops $L_k$:
$A' x + B' y  = k/D \mod 1$ with $k \in \{0,1,\dots,D-1\}$. 

Consider $R^{nk}_{\varphi,\vartheta}(0,0)$.  
Since $nD(A'\phi + B'\theta) = nD(\frac{C}{D}) = 0 \mod 1$
these points are contained in the torus-loop $L_0$.  As in the case
$\Delta = 3$ the action of $R^{nD}_{\varphi,\vartheta}$ on the circle
$L_0$ is an irrational rotation, therefore the orbit of $(0,0)$ is
dense in $L_0$.

Suppose that $R^{n}_{\varphi,\vartheta}(0,0) \in L_k$.  Then
we claim that $R^{n+1}_{\varphi,\vartheta}(0,0) \in L_{k'}$ with
$k' := (k + C) \mod D$.
To see this note that $n(A'\varphi + B'\vartheta) = \frac{k}{D}$ implies
that 
$$(n+1)(A'\varphi + B'\vartheta)  = \frac{k}{D} + 
(A'\varphi + B'\vartheta) = \frac{(k + C) \mod D}{D} = \frac{k'}{D}.$$  
Since $\gcd(C,D) =1 $ the map $R^{n}_{\varphi,\vartheta}$ permutes the
circles $L_k$ cyclically.  Thus the orbit of $(0,0)$ is dense
in $\L := \cup_{k=0}^{D-1} L_k$. 

We claim that if $A,B,C$ are sufficiently large that then 
$\L$ is $\e$-dense in $\T^2$.  There are five subcases.

i) Suppose $B' = 0$.  The circles $L_k$ are vertical, and are given
by the formula $x = k/|A|$.  If $|A| > 1/\e$ then these circles are 
$\e$-dense.

ii) Suppose $A' = 0$.  The circles $L_k$ are horizontal, and the
proof is similar to case i).

iii) If $B' \ne 0$ and $|A'| > 1/\e$ then each 
then each of the circles $L_k$ is $\e$-dense and the lemma follows.

iv) The case $A' \ne 0$ and $|B'| > 1/\e$ is similar
to case iii).

v) If $D$ is large we have $D$ parallel and equidistant circles, thus
they are at least $1/D$ dense.

The only cases which do not fall into these five subcases are
when $|A'|,|B'|$ and $D$ are all smaller than $1/\e$, 
and thus $|A| = |A'D| < 1/\e^2$  and $|B| = |B'D| < 1/\e^2$.
Furthermore $|C| \le |A| + |B| < 2/\e^2$.  
Thus the exception set of $\theta$ is finite.
\end{proof}


\begin{proofof}{Theorem \ref{thm2}} Throughout the proof $p$ and $q$
stand for points in $\abx$.
We call a path connecting $p$ to $q$ a $VH$-path if it is the union  of
a finite number of vertical and horizontal segments. We will call a 
point in between a vertical and horizontal segment a turning point. 
We claim that in a prelattice surface every $p$ and $q$ can be connected by
a finite length $VH$-path.  We give an algorithm to produce simple 
$VH$-paths between
two points.
First each $q$ in $\CC_v(p)$ (resp.~$\CC_h(p)$) can be connect to $p$ by a
$VH$-path with one turning point. Next each $q' \in \CC_v(q)$ or 
$\CC_h(q)$ can be connected to $p$ with a $VH$-path with two turning points,
(see Figure 2) etc. 
\medskip\\

\begin{figure}[h]
\psfrag*{t}{{\footnotesize{turning point}}}
\psfrag*{p}{{\footnotesize{$p$}}}
\psfrag*{q}{{\footnotesize{$q$}}}
\centerline{\psfig{file=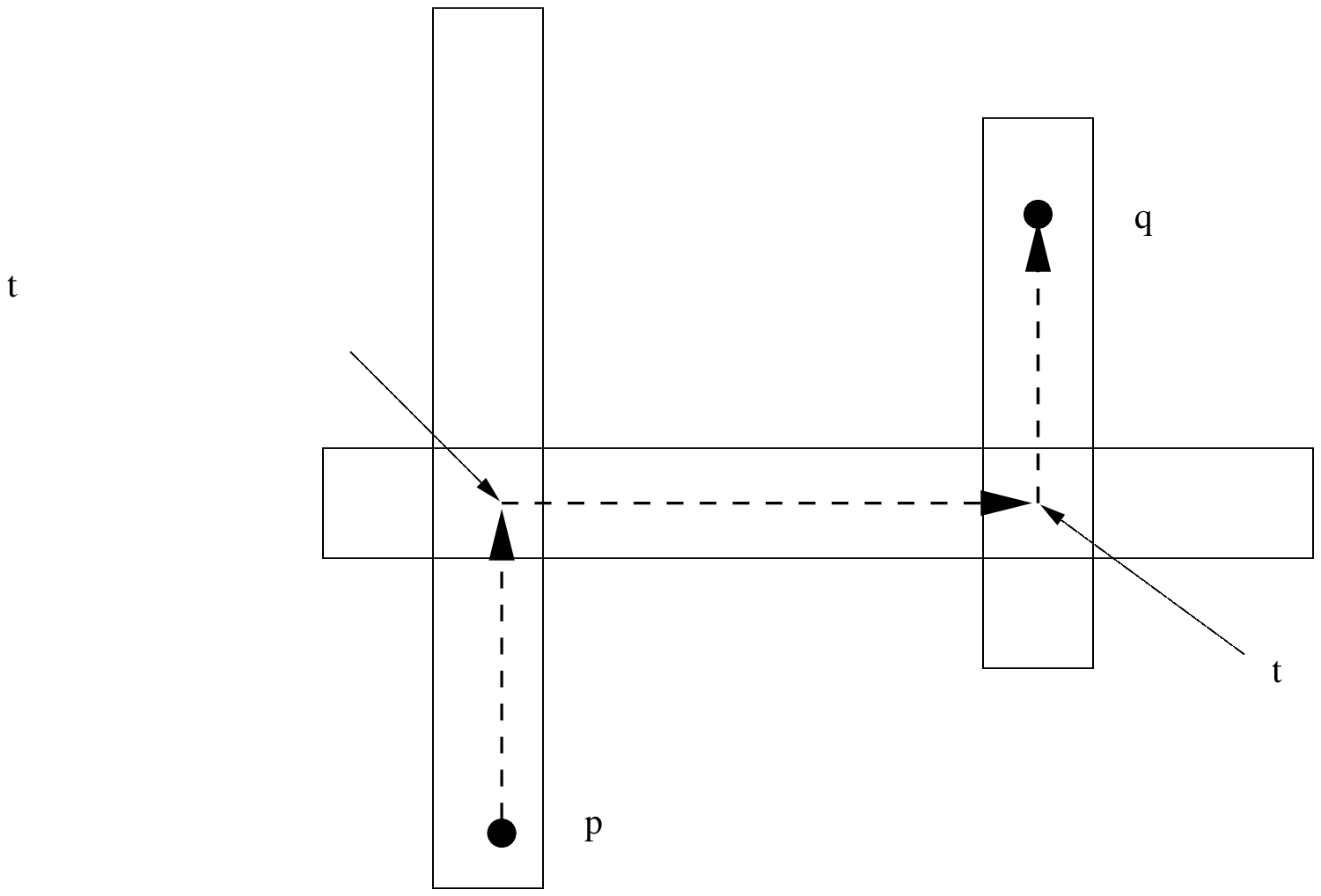,height=44mm}}
\caption{A VH-path.}
\refstepcounter{figure}\label{VH}
\end{figure}

Since the surface $\abx$ is connected 
the claim follows by exhaustion. From here on
a $VH$-path will always be one of the simple paths produced by this
algorithm.  
Fix $p$ and $q$.  Let $r_i$ be the sequence of turning points 
of
a $VH$-path connecting $p$ to $q$ and let $\Rec_i := \Rec(r_i)$.
Clearly the sequence $\Rec_i$ depends only on $\Rec(p)$ and $\Rec(q)$ rather
than on $p$ and $q$.  We call this sequence of rectangles the 
combinatorics of the rectangles $\Rec(p)$ and $\Rec(q)$.

Our strategy is as follows.  We fix a point $p$ and a 
rectangle $\Rec = \Rec(p)$.
We are interested in the $q \in \Rec_N$ which are illumined by $p$.
We consider powers of $T_h$ and $T_v$ which keep $p$ in $\Rec(p)$ and
try to move $q$ closer to $p$.  More precisely fix
the combinatorics of $\Rec(p)$ and $\Rec(q)$.  By induction it suffices to
show that  all but countably many 
points in the rectangle $\Rec_N$, can be moved into rectangle $\Rec_{N-1}$
by a power of $T_h$ or $T_v$ which keeps $p$ in $\Rec(p)$.

To prove this we apply Lemma \ref{Kronecker}.  Suppose that $\Rec_N$ 
and $\Rec_{N-1}$ are connected by a vertical path (i.e.~they are in the 
same vertical cylinder $\CC$ of height $h$ and width $w$). Remember that
if $\CC$ is a vertical cylinder then on each vertical segment
$T_v$ acts as a rotation. We consider the vertical
cylinders $\CC_v(p) \times \CC_v(q)$ and in this set we consider the
product $\T^2 = \S^1 \times \S^1$ consisting of the vertical
closed curve through $p$ times the vertical closed curve through 
$q$.\footnote{Note that as $q$ varies in $\Rec(q)$ these tori are isometric
copies of each other, thus we can apply the abstract Lemma \ref{Kronecker} 
as if they are all the same.}  
The action of $T_v \times T_v$ is a translation on this torus.
The angles of rotation are given by \eqref{vangle}.
Thus we can apply Lemma \ref{Kronecker} with
normalized coordinates $c := x/h$, $d := y/w$ and $\e$ chosen
such that $(c-\e,c+\e)$ corresponds (via the normalization) to 
the intersection the vertical closed curve through $p$ and 
the rectangle $\Rec(p)$.
Furthermore
$(e-\e,e+\e)$ corresponds (via the normalization) to 
the intersection the horizontal closed curve through $p$ and 
the rectangle $\Rec(q)$.
We conclude that except for an exceptional set consisting of a
finite union of vertical segments 
of $\Rec_N$ all other points in $\Rec_N$
can be moved to $\Rec_{N-1}$ while the image of $p$ stays in $\Rec(p)$.
We call this union of vertical segments $\CB_v(p)$.

Next we must apply now the horizontal twist.
We apply Lemma \ref{Kronecker} to the torus
$\T^2 = \S^1 \times \S^1$ in $\CC_h(p) \times \CC_h(q)$ with $\e$ chosen
so that $(c-\e,c+\e) \subset \Rec(p)$ and $(e-\e,e+\e) \subset \Rec_N$.
We conclude that  except for an exceptional set consisting of a
finite of horizontal segments of $\Rec_N$ all other points in $\Rec_N$
return to  $\Rec_{N}$ while the image of $p$ stays in $\Rec(p)$.
We call this union of horizontal segments $\CB_h(p)$.

For each $k \ge 1$ consider the set 
$X_k := \{q \in \Rec_N: T^k_h(q) \in \Rec_N\}$.
Thus  $\cup_{k \ge 1} X_k = \Rec_N \setminus B_h(p)$. The collection of
points $\CB_h(p) \cap \CB_v(p)$ is finite.  
These points are exceptional points.

Applying the previous argument with the vertical twist 
yields that except of an exceptional set 
consisting of a finite  union of vertical 
segments $\CB_v(T^k_h(p))$ in $X_k \subset \Rec_N$
all other points in $X_k$ can be moved to $\Rec_{N-1}$ while the image of 
$T^k_h(p)$ stays in $\Rec(p)$.  Applying $T^{-k}_h$ to $\CB_v(T^k_h(p))$ 
yields a finite collection of segments of slope $-1/ak$ in $\Rec_{N}$. 
Since the sets  $\CB_v(q)$ and $T^{-k}_h(\CB_v(T^k_h(p)))$ consist of a finite
union of segments which intersect transversely
the intersection
$\CB_v(q) \cap T^{-k}_h(\CB_v(T^k_h(p)))$ is a finite set and the points in
$X_k$ which are not in this finite set can be moved 
to $\Rec_{N-1}$ while the image of $p$ stays in $\Rec(p)$.
We conclude the theorem by taking the union over $k$.
\end{proofof}

\section{Illumination on Veech surfaces}

\begin{proofof}{Theorem~\ref{Schmoll}}
To prove Theorem~\ref{Schmoll} note that a
Veech surface $\abx$ is necessarily a prelattice surface, thus we can
apply Theorem \ref{thm1}.

Let us assume that $X$ is a non arithmetic Veech surface.
Given $p \in X$ we consider any convex open open neighborhood  $\CO_p$ 
of $p$ and the {\em exceptional set}
$$\sE:= \big \{(z_h,z_v) \in X \times X: 
\aff \abx \diact (z_h,z_v) \not \in \CO_p \times \CO_p  \big \}.  $$
The set $\sE$ is clearly closed and $\aff \abx$ invariant.
Our results on invariant subspaces say, that 
$\sE$ is a finite union of orbit closures 
$\overline{\aff\abx \diact (z_1,z_2)}$, 
of complex dimension less than $2$.
For $p \in \CO_p$, the set $\sE \cap (\{p\} \times X)$ is algebraic. Additionally
   by Theorem \ref{thm1}, this set is at most countable, thus
 it is finite.

 Let us assume now that $X$ is an arithmetic surface.
 
 We first suppose that $p$ is a periodic point which means that both coordinates of $p$ are rational. In this case, we refine the proof of theorem \ref{thm1} to prove that $p$ illumines every point $q$ except a finite number.
In the sequel, we use the notations of the proof of Theorem \ref{thm1}.  Without loss of generality, we  assume that $p$ is fixed by  horizontal and  vertical parabolic elements $T_h$ and $T_v$.
Let $\Rec(q)$ be the rectangle defined by the horizontal and vertical cylinder decompositions of $X$ containing $q$.  We consider a VH-path from $\Rec(p)$ to $\Rec(q) = \Rec_N$. By induction, we show that all but {\it finitely many} points in the rectangle $\Rec(q)$ can be moved to $\Rec_{N-1}$.
Suppose that $\Rec_N$ and $\Rec_{N-1}$ are in the same vertical cylinder. 
As in the proof of Theorem \ref{thm1}, applying $T_v$, we can move to $\Rec_{N-1}$ every point of $\Rec_N$, except a finite number of segments denoted by $\CB_v(p)$. The vertical coordinates of the elements of $\CB_v(p)$ are rational numbers with small denominators. If $q$ belongs to $\CB_v(p)$, we apply the horizontal twist $T_h$.  If the $x$ coordinate of $q$ is irrational, the orbit of $T_h$ is dense in the horizontal circle containing $q$. If this coordinate is rational with a denominator large enough, it is $\varepsilon$ dense. If $\varepsilon$ is chosen small enough,
 the orbit of $q$ under $T_h$ intersects $\Rec_n \setminus \CB_v(p)$.
 Therefore, except for a finite number of points in $\CB_v(p)$, this orbit intersects  $\Rec_n \setminus \CB_v(p)$. Thus, applying $T_v$, we move $q$ to $R_{N-1}$.
 
 If $p$ is not a periodic point, we consider the exceptional set $$\sE:= \big \{(z_h,z_v) \in X \times X: 
\aff \abx \diact (z_h,z_v) \not \in \CO_p \times \CO_p  \big \},$$
\noindent where $\CO_p$  is a convex open neighborhood of  p. As $p$ is not periodic, the fiber 
$\sE_p = \sE \cap (\{p\} \times X)$ does not contain any periodic point. Moreover, the intersection of $\sE_p$ with an invariant two dimensional orbit closure  is finite. Therefore if the set $\sE_p$ is infinite, $\sE$ contains an infinite number of  invariant two dimensional surface. Let $\pi$ be the projection from $X$ to $\T^2$.  The projection of $\sE$ is not $\slz$ invariant. But, as $\slv{X,\omega}$ is a finite index subgroup of $\slz$, there exists a set $Y$, closed and $\slz$ invariant, contained in $\T^4$ which is a {\it finite} union of images of $\pi(\sE)$.
The set $Y$  contains infinitely many two dimensional torii. By proposition \ref{t4inv}, it is equal to $\T^4$. This leads to a contradiction because Theorem \ref{thm1} implies that the fiber over $\pi(p)$ is at most countable. Thus, the set   $\sE_p$ is finite which completes the proof of theorem \ref{Schmoll}.

\end{proofof}
\medskip\\
\begin{proof}[{\bf Proof of Theorem \ref{nonilluconf}}]
We first prove the
\begin{lemma}\label{self-illumination} Let $X$ be a Veech surface then every point $p \in $X  which is non periodic 
is self-illuminating.
\end{lemma}
\begin{proof}

Let $p$ be a non periodic point. We take a periodic direction, say the 
horizontal one, and decompose $X$ into maximal cylinders. 
Now we fix an open maximal cylinder $\CC$, 
in this decomposition. 
Then the horizontal leaf through any $m \in \CC$ is regular and thus $m$ illuminates itself. 
If $p \in X$ is not periodic its $\aff \abx$ orbit 
is dense, hence there is a point, say $\psi(p)$  
in the orbit of $p$ which is in $\CC$ 
and therefore the horizontal leaf through $\psi(p)$ is 
regular. Thus  $p$ 
is self-illuminating.   
\end{proof}

By our results on invariant subspaces in arithmetic surfaces  
$\abx$ cannot be arithmetic.
Without restrictions we assume that $\abx$ is reduced. 
Recall that a point $(p,q)\in X^2$,  represents an illuminable configuration, 
i.e. there is a saddle connection 
between $p$ and $q$, whenever its $\slv{X,\omega}$ orbit intersects 
a convex set $\CO \times \CO \subset X^2$. In the sequel, we will assume that $(p,q)$ is not a periodic point.

    Now our classification says either
\begin{itemize}
\item  $(p,q) \in X^2$ is generic and $D_+\subset \overline{\slv{X,\omega}\diact (p,q)}$, or  
\item  $(p,q)$ is part of an invariant surface $(S,\omega_S)$.
\end{itemize}

All points in $X^2 \backslash D_+$ whose $\slv{X,\omega}$ orbit closure 
contain the diagonal $D_+$ represent illuminable configurations. This
is obvious if $p\neq q$ and a direct consequence of 
Lemma \ref{self-illumination} if $p= q.$

If $D_{\pm}$ are the only $\slv{X,\omega}$ invariant surfaces 
of slope $\pm 1$, then all invariant surfaces in $X^2$ intersect $D_+$. 
For horizontal or vertical subspaces this is clear.  
Denote the set of fix points 
under the (unique) affine involution by  
$Fix_{\pm}:= D_- \cap D_+$.


Because $Fix_{\pm}$ is an intersection of 
invariant surfaces, all its points are periodic points. 
Moreover the diagonal action of $\abx$ on $X^2$ restricted to 
invariant surfaces is ergodic too. This implies immediately     
that all nonperiodic points on parallel invariant surfaces 
are illuminable configurations. 
The only question remaining is what happens on $D_-$. 
The answer is: 
\begin{itemize}
\item if all points in $Fix_{\pm}$ 
are cone points, every nonperiodic point $(p,q)=(p,\phi(p))\in D_-$ 
defines a non illuminable pairing  
\item if $Fix_{\pm}$ contains one regular point all non periodic 
points on $D_-$ define illuminable configurations on $X$. 
\end{itemize}

The last statement follows again from ergodicity of 
the $\slv{X,\omega}$ action restricted to $D_-$ and from the

\begin{lemma}\label{fix} Let $\phi \in \aff \abx$ be the affine involution. The set
$Fix_{\pm}$  contains regular points if and only if 
there is  $p \in X \backslash Z(\omega)$ and a saddle connection from $p$ to $\phi(p)$.
\end{lemma}

\begin{proof}
Given $p \in X \backslash Z(\omega)$ 
with $\phi(p)\neq p $, every saddle 
connection from $p$ to $\phi(p)$ contains 
a fix point of $\phi$, since $\phi(s)=-s$ 
($-s$ is $s$ with the opposite orientation).
This follows because $\phi$ exchanges $p$ 
and $\phi(p)$,  while the angle of any outgoing 
or incoming saddle connection is turned by 
$\pi$. 

The converse is obvious.
\end{proof}

Finally note that if $X$ is not reduced $\pi^{-1}(Fix_{\pm})\subset X$ 
might only contain cone points, while $Fix_{\pm}$ does not.   
Here $\pi: X \rightarrow X_{red}$ is the canonical cover. 
\end{proof}

\begin{proof}[{\bf Proof of Corollary \ref{cornonilluconf}}]
By Corollary \ref{knownexample} the Veech surfaces in view 
admit only the diagonal and off-diagonal as 
invariant subspaces. Since there are always regular Weierstrass points on genus $2$ surfaces with one cone point, 
one only needs to consider pairs of periodic points.  
Moreover by Möllers 
result \cite{mm} we know that the only 
periodic points on non-arithmetic Veech surfaces of genus $2$ 
are the Weierstrass points.
\end{proof}
{\bf General non-illumination or blocking configurations.}
Let  $\abx$ a Veech surface with an involution 
$\phi \in \aff \abx$,  then   
\[ D_{\phi}=\{ (z,\phi(z))\in X^2: z \in X\} \]
is an invariant subset of slope $-1$ in $X^2$. 
By the argument in the Proof of Theorem \ref{nonilluconf}  
$Fix_{\phi}$ provides a {\em blocking set} for $SC(z,\phi(z))$.   
The set $Fix_{\phi}$ consists of finitely 
many periodic points because $Fix_{\phi}=D_{\phi} \cap D_+$, 
hence we obtain a {\em finite blocking property} 
for all non-periodic pairs contained in $D_{\phi}$. 
Because $D_{\phi} \cong X$ as translation surfaces, {\em generic} 
is equivalent to {\em not periodic} with respect to 
the action of the Veech-group.
\medskip\\
{\bf Example.}
We again look at the staircase, but 
now we take only staircases with an even number $2n$ 
of steps. Denote by $(Y_n,\omega_n)$ the staircase 
with $2n$ steps, then we have the following: 
\begin{itemize}
\item The absolute periods define a covering 
$\pi_n:(Y_n,\omega_n) \rightarrow (\R^2/2\Z^2,dz)$ 
branched over the points $\Z^2/2\Z^2$  
\item $(Y_n,\omega_n)$ has $4$ cone points of order 
$n$, or equivalently $\omega_n$ has $4$ zeros of 
order $n-1$
\item $g(Y_n)=2n-1$ 
\item $\aut (Y_n,\omega_n)\cong \Z/4n\Z $
\item $\aut_{\pi_n}(Y_n,\omega_n):=\{\phi \in \aut (Y_n,\omega_n): \pi_n\circ\phi=\pi_n\}
\cong \Z/n\Z $
\item $\slv{Y_n,\omega_n} \cong \Gamma_2$ unless $n=1$, but then
\item $(Y_1,\omega_1) \cong (\T^2,dz)$ and $\slv{Y_1,\omega_1} \cong \slz$ 
\end{itemize}

\begin{lemma}
For any point $z \in Y_n\backslash Z(\omega_n)$ there are $n$ points, 
given by $\pi^{-1}_n(-\pi_n(z))$, which are not illuminable 
from $z$.  
\end{lemma}
\begin{proof}
Assume there is a connecting geodesic, say $s$ 
between $z$ and $z_{-} \in \pi^{-1}_n(-\pi_n(z))$. 
Then $\pi_n(s)$ connects $\pi_n(z)$ and $\pi_n(z_-)=-\pi_n(z)$ 
and therefore contains a Weierstrass point  
of $\R^2/2\Z^2$. That means $s$ intersects one of the 
preimages of Weierstrass points of $\R^2/2\Z^2$, 
but these are all cone points.
\end{proof}
Consequently, there are surfaces $\abx$ such any point $z \in X$,  
which is not a cone point, admits  
several points not illuminable from $z$. In particular the total 
number of points which are not illuminable (from a given point) 
has no universal bound. 
Since the Veech-group of this example has not been used  
we drop it and present a more general class of surfaces 
with non-illumination pairs. 
\medskip\\
{\bf Non-Veech example with non-illumination points.} 
Take any surface represented by an $L$ shaped figure 
and construct the following double cover.
\begin{figure}[h]
\epsfxsize=9truecm
\centerline{\epsfbox{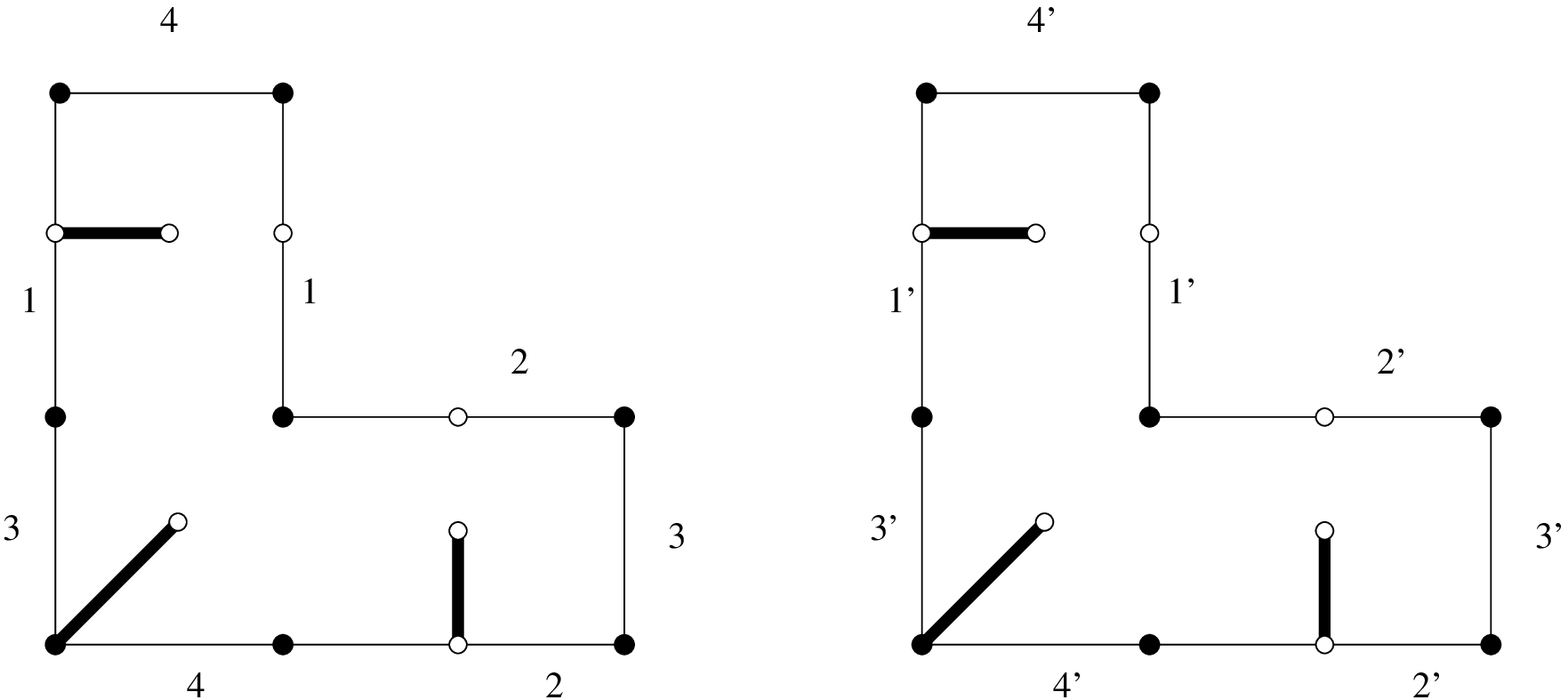}}
\caption{Double cover of $L$ branched over Weierstrass points.}
\refstepcounter{figure}\label{lcover}
\end{figure}
\vspace*{1mm}\\
To obtain this cover, $L_C$, 
cut each copy of the $L$ along the thick lines, 
and identify the two parallel slits so obtained cross wise. 
No matter what the width-height relations of $L$ are, 
$L_C$ is a double cover of $L$ with cone points above  
all Weierstrass points of $L$. Recall that 
each $L$-shaped surface admits an involution.  

By our previous arguments all straight 
geodesics connecting pairs of regular 
points $p,q \in L_C$, which are preimages of pairs $a,\phi(a) \in L$ 
in involution, are crossing a cone point. 
\medskip\\
{\bf Polygonal Billiards.} 
There is a difficulty to construct examples of polygonal 
billiard along these lines: the unfolding usually interferes 
with the covering construction. More precisely: 
given an $n$-gon $G$ and an unfolding $UG$ with 
unfolding group $H$ 
and suppose there is a covering map $\pi: UG \rightarrow \R^2/2\Z^2$, 
then given a point $p \in UG$   
the set $F_p :=\pi^{-1}(-\pi (p))$ is not always invariant under 
the action of $H$. At least we were not able to  
find an example for which $H \cdot F_p= F_p$.  
Connecting $p$ with a point $q \in H \cdot F_p \backslash F_p$ 
one can  fold the obtained geodesic segment to a saddle connection    
connecting $p$ with a point $F_p \cap G$.  
This destroys examples like the swiss cross 
tiled by five squares.    
\medskip\\

\noindent {\bf Other applications of invariant subspaces.}

The initial purpose of looking at $\slv{X,\omega}$ invariant 
subspaces $S \subset X^2$ of dimension $2$ was using  
their nice structure as a translation surface to 
write down asymptotic constants in terms of the translation 
geometry of the invariant subspace $S$, like in \cite{EsMaWi} 
or in \cite{s2,s3}. In particular we can look at modular fibers 
$\sF$ of surfaces parametrizing branched covers of $X$, 
such that there is a covering map $\pi: \sF \rightarrow S$. 
Note that in such a case $\sF$ is a lattice surface $(\sF,\omega_{\sF})$ 
and
\[\slv{\sF,\omega_{\sF}}\cong\slv{S,\omega_{S}}\cong\slv{X,\omega}.\] 
We postpone results regarding quadratic asymptotic 
constants in connection with invariant subspaces of $X^n$ 
to the forthcoming paper \cite{s5}. 
 Theorem \ref{veech-groups} together 
with the other results on invariant subspaces 
and the asymptotic formulae described in \cite{s2,s3} 
allow to evaluate asymptotic constants without  
using the Siegel-Veech formula directly.


\section{Further remarks and open problems}

As pointed out to us by Thierry Monteil, for a ``generic''
translation surface every point illumines every point.
For example, for hyperelliptic strata, Veech has
shown that for almost all parameter choices one can choose
a fundamental domain which is convex \cite{v2}, and thus every point
illumines every point.  A similar but more technical argument
works in other strata.

Alex Eskin remarked to us that our type of theorem will hold
in any case to which we can generalize Ratner's machinery.

All known examples of nonillumination in polygons use the
same mechanism which leads us to state the following 
conjecture.

Conjecture 1:  In any rational polygon $P$, any point $p$
illumines all but finitely many points $q \in P$.
\medskip\\
Regarding our results, in particular that non-illuminable pairs 
are already pairs of periodic points, one might suspect 
that on a general translation surface, i.e. not a branched cover of 
a Veech surface, there are no nonilluminable 
configurations, unless they are (indirectly) caused by 
fixed points of an involution.  
\medskip\\
The convention we have taken is that orbits which arrive at corners
of a polygon (singularities of a translation surface) stop.  There
are always (at most) two continuations possible, which we call
``singular orbits or geodesics''.

Conjecture 2: Consider a translation surface $\abx$ and $p,q \in X$
such that $p$ does not illumine $q$.  Then there is a ``singular geodesic''
connecting $p$ to $q$.

It is natural to study the points which are illumined at time $t$,
i.e.~the range of the map $\exp_p(t)$.
This leads to the following questions.

Question 1:  Consider a translation surface $\abx$. For a point $p \in X$ 
consider the exponential map $\exp_p(t)$. For which $p$ is the
set $\exp_p(t)$ asymptotically dense? asymptotically well distributed?

Question 2: What can be said for irrational polygons?

\end{document}